\documentclass[12pt,a4paper,leqno]{article}
\usepackage{a4wide}
\usepackage{latexsym}
\usepackage{amsmath}
\usepackage{amssymb}
\newtheorem{definition}{Definition}
\newtheorem{lemma}{Lemma}
\newtheorem{proposition}{Proposition}
\newtheorem{theorem}{Theorem}
\newtheorem{corollary}{Corollary}
\newenvironment{proof}{\medskip\par\noindent{\bf Proof}}{\hfill $\Box$
\medskip\par}
\begin{document}
\title{On the generalized Riemann-Hilbert problem with irregular
singularities}

\author{A. A. Bolibruch\addtocounter{footnote}{4}\footnote{Steklov Mathematical Institute, Gubkina str. 8, 117966 Moscow, Russia . During the preparation of the paper, the late
Andrey Bolibruch was an invited professor   at  Institut de Recherche Math\'ematique Avanc\'ee, Universit\'e Louis Pasteur et CNRS,
7 rue Ren\'e Descartes, 67084 Strasbourg Cedex, France.},\ \    S. Malek\footnote{Universit\'e de Lille 1, UFR de math\'ematiques, 59655 Villeneuve
d'Ascq Cedex France, email: Stephane.Malek@math.univ-lille1.fr. During the preparation of this paper, the author was supported by
a Marie Curie Fellowship of the European Community programme ``Improving Human Research Potential'' under contract number
HPMF-CT-2002-01818.} \ and \ C. Mitschi\footnote{Institut de Recherche Math\'ematique Avanc\'ee, Universit\'e Louis Pasteur et CNRS,
7 rue Ren\'e Descartes, 67084 Strasbourg Cedex, France, email: mitschi@math.u-strasbg.fr}}
\date{}
\maketitle

\def\cf{{\it{cf.}}}
\def\ff{{\mathcal F}}
\def\aa{{\mathcal A}}
\def\gg{{\mathcal G}}
\def\pp{{\mathbb{P}^1(\mathbb{C})}}
\def\gl{{\mathrm{GL}}}
\def\cc{\mathbb C}
\def\zz{{\mathbb Z}}
\def\qq{{\mathbb Q}}
\def\ll{{\mathcal L}}
\def\cn{{\mathcal N}}
\def\ct{{\mathcal T}}
\def\cd{{\mathcal D}}
\def\nn{\nabla}
\def\dd{\mathrm {deg}}
\def\diag{\mathrm {diag}}
\def\rr{\mathrm {rank}}
\def\mm{{\mathcal M}}
\def\ss{{\mathcal S}}
\def\ee{{\mathcal E}}
\def\tt{{\check{T}}}
\def\sb{\overline{s}}
\def\mb{\underline{M}}
\def\rb{\overline{r}}
\def\ee{\boldsymbol{\mathcal{E}}}

\renewcommand{\abstractname}{Abstract}
\begin{abstract}
In this article we study the  generalized Riemann-Hilbert problem, which extends the classical Riemann-Hilbert problem to the
case of irregular singularities. The problem is stated in terms of  generalized monodromy data which  include
the monodromy
representation, Stokes matrices and the true Poincar\'e rank at each singular point. We give sufficient conditions for the existence of a
linear
differential system with such data. These conditions are in particular fulfilled when the monodromy representation is irreducible, as in
the classical
case. We solve the problem almost completely in dimension two and three. Our results have  applications in  differential Galois theory. We give sufficient conditions for a given linear algebraic group $G$ to be the differential Galois group over $\cc(z)$ of a linear differential system with a minimum number of singularities, all fuchsian but one, at which the system has a minimal Poincar\'e rank.\end{abstract}

\noindent
There are many approaches to differential equations. One can focus on the existence and behaviour of  the solutions, or on algebraic
properties of their
symmetries. One may also ask for the existence of  differential equations that satisfy specific inverse problems such as the
Riemann-Hilbert problem,
the Birkhoff standard form problem or the inverse problem in differential Galois theory. This article is an attempt to relate the three
problems through
the statement and solutions of the generalized Riemann-Hilbert problem.

\medskip
\noindent 
The classical Riemann-Hilbert problem asks for conditions under which a given representation 
$$ \chi:\pi_1(\mathbb{P}^1(\mathbb{C})\setminus\mathcal{D},z_0)\longrightarrow {\mathrm {GL}}(p;\mathbb C) $$
of the fundamental group of the Riemann sphere ${\mathbb{P}}^1({\mathbb{C}})$ punctured at each point of a finite subset $\cal D$ not
containing $z_0$,
can be realized as the monodromy representation of a linear differential system with fuchsian singularities only, all in $\cal D$. Let
us recall that
a point $a\in \cal D$ is a {\it fuchsian} singularity of a linear differential system ${dy}/{dz}=B(z)y$, where $ B$ is an  $n\times n$
matrix
with coefficients in $\cc (z)$,  if $a$ is a simple pole of $B$ ({\it modulo} a M\"obius transformation if $a=\infty$).
This problem is still open, although  important results of A. Bolibruch (\cite{Bol1}, \cite{Bol2}, \cite{Bol3}, \cite{Bol4}) have reduced it
considerably. Several authors have given sufficient conditions either to solve this problem or to construct counterexamples.
A. Bolibruch \cite{Bol1} and V. Kostov \cite{Ko} have shown independlently that the irreducibility of the representation $\chi$ is
a sufficient condition. In dimension two the problem always has a solution (\cf \cite{AB}) and in dimension three and four it has been completely elucidated (\cite{AB}, \cite{Bol2}, \cite{Glad}).
The Riemann-Hilbert problem is related to problems in many  areas of mathematical physics and has become a trend of research over the last twenty years. There is extensive literature available on the subject, in particular on Painlev\'e equations and isomonodromic deformations. For recent results in this field we refer to  
\cite{Boa}, \cite{Bol6}, \cite{DM}, \cite{Hi}, \cite{It}, \cite{Iso}, \cite{KiKo}, \cite{NO}.
 
\medskip
\noindent 
Closely related to the Riemann-Hilbert problem, the Birkhoff inverse problem asks the following. Consider a differential system $z dy/dz=A(z)y$
where the matrix $A(z)=z^r\sum_{n=0}^{\infty}A_nz^{-n}$ is meromorphic at infinity. Does there exist a differential system
$z dy/dz=B(z)y$,  where $B(z)$ is a  polynomial coefficient matrix, meromorphically equivalent to the given system and
with a Poincar\'e rank at infinity not greater than the original one?
In dimension two and three, the problem is known to have a positive answer, see \cite{Ba2}, \cite{JLP1},
but for in higher dimension, although many sufficient conditions have been given, see \cite{Bol3}, \cite{BB}, \cite{Sa},
the problem remains open in general. The differential systems in Birkhoff standard 
form appear in complex algebraic geometry in the study of particular Frobenius manifolds, see \cite{Sa} and
references therein.

\medskip
\noindent
In the present paper we extend both the Riemann-Hilbert problem and the Birkhoff standard form problem to the case of an arbitrary number of irregular singularities. We define generalized monodromy data, consisting of the monodromy representation with respect to  prescribed singularities and  of further prescribed local data at each singularity. These data include the Poincar\'e rank and  Stokes data. The generalized Riemann-Hilbert problem is the following: {\it Let singular points and generalized monodromy data be given in which all Poincar\'e ranks are minimal. Construct a system  on $\pp$  with these data}. We give sufficient conditions to solve this inverse problem and we show that they are in general fulfilled in dimension two and three.

\medskip
\noindent 
We conclude the paper with applications to  differential Galois theory, where we under suitable assumptions solve the inverse problem with a better control of the singularities.  The global inverse problem in differential Galois theory over ${\mathbb{P}}^1({\mathbb{C}})$, that is, over the differential field $K=\mathbb{C}(z)$, asks for the existence of a differential system $dy/dz=B(z)y$ with coefficients in $\mathbb{C}(z)$ and with a given linear algebraic group (over $\mathbb{C}$) as its differential Galois group over  $\mathbb{C}(z)$. It  always has a solution. This was first proved by  M. and C. Tretkoff \cite{TT}, using a weak solution of the Riemann-Hilbert problem. Other proofs were given since,  either analytic  (\cite{R4}, \cite{R5}, \cite{R6}), or algebraic over a general  field of constants (\cite{JH}). An algebraic and constructive proof  was given  in \cite{MS1} for connected groups, and in \cite{MS3}, \cite{CMS} for large classes of non-connected groups. In the present paper we focus on the number and on  the Poincar\'e rank of the singularities of  a differential system with a given  Galois group, and we show that under suitable conditions both  are minimal. 

\medskip
\noindent 
The paper is organized as follows.

\medskip
\noindent In section 1  we define  generalized monodromy data attached to a linear differential system over $\pp$.

\smallskip
\noindent In section 2  we state the generalized Riemann-Hilbert problem and we show that it has a solution if a certain family of vector bundles with connections contains a holomorphically trivial  bundle.

\smallskip
\noindent In section 3 we give further sufficient conditions, in terms of the stability of a certain bundle,  to solve the problem. These conditions are in particular fulfilled when the monodromy representation is irreducible and the data at one of the singularities are unramified. If all data are fuchsian, we recover the irreducibility condition of Bolibruch and Kostov.

\smallskip
\noindent In  section 4 we look for the existence of possibly reduced systems with given generalized monodromy data, when the prescribed monodromy data are ``non generic''. This should lead to a reduction of the problem to an equivalent problem in lower dimension. 

\smallskip
\noindent In section 5 the results for reduced systems enable us to  solve the generalized Riemann-Hilbert problem  completely in dimension two and three, assuming that not all the singularities are irregular with ramification.

\smallskip
\noindent In  section 6 we apply our previous results to  the inverse problem of differential Galois theory, which under suitable conditions can be solved with a minimal number of singularities and a minimal  Poincar\'e rank at these.

\section{Generalized Monodromy Data}
Consider a system
\begin{equation}   \label{s}
\frac{dy}{dz}=B(z)y
\end{equation}
of $p$ linear differential equations with rational coefficients on
the Riemann sphere $\mathbb{P}^{1}(\mathbb{C})$. Let $\mathcal{D} = \{ a_{1},\ldots,a_{n} \}$ be the set of singular points of (\ref{s}), consisting of 
the poles of the matrix function
$B(z)$ and of a possible singular point at infinity (if the system obtained from (\ref{s}) {\it via} $z={1}/{u}$ has a singular point at the origin).

\noindent
Consider the matrix differential form $\omega=B(z)dz$. In what follows we will rather write (\ref{s}) in its invariant form
\begin{equation}   \label{sinv}
dy=\omega y.
\end{equation}
in terms of which   $\mathcal{D}$ is  a singular divisor of $\omega$.

\medskip
\noindent
To any system (\ref{sinv}) there correspond what we will call {\it generalized
monodromy data}, which we define below.

\subsection{The monodromy representation}
Let $Y$ denote a fundamental solution of (\ref{sinv}), holomorphic in a neighbourhood of a given non-singular point $z_0\in\mathbb{C}$. Analytic continuation of $Y$ along a loop $\gamma$ in $\mathbb{P}^{1}(\mathbb{C})\setminus\mathcal{D} $ yields a new fundamental solution $\gamma^*(Y)=YG_{\gamma}$ for some matrix $G\in\mathrm{GL}(p,\mathbb{C})$. This defines the {\it monodromy representation}    
\begin{equation}   \label{rep}
\chi:\pi_1(\mathbb{P}^1(\mathbb{C})\setminus\mathcal{D};z_0)\longrightarrow
{\mathrm {GL}}(p,\mathbb C)
\end{equation}
of the system, with respect to $Y$. 
Since the fundamental group of $
\mathbb{P}^1(\mathbb{C})\setminus\cal{D}$ is generated by the homotopy
classes of all elementary loops $\gamma_i$, where $\gamma_i$, $i=1,\ldots,n$, encloses  the only singular point $a_i$, the monodromy representation of (\ref{sinv}) is
defined by the local monodromy matrices $G_i$   corresponding to these
loops. These matrices satisfy {\it a priori}  the only relation
$G_1\cdot\ldots\cdot G_n=I$.

\subsection{The Poincar\'e rank}
Let $a\in\mathcal{D}$ be a given singular point $a_i$ of (\ref{sinv}) and $G$ the corresponding monodromy matrix $G_i$. In the neighbourhood of $a$ the coefficient matrix of (\ref{s}) can be expanded as follows

\begin{equation}   \label{rP}
B(z)=\frac{B_{-r-1}}{(z-a)^{r+1}}+\ldots+\frac{B_{-1}}{z-a}+B_0+\sum_{i=1}^\infty {B_{i}}{(z-a)^{i}}
\end{equation}
where $B_{-r-1}\neq 0$.
\begin{definition} The {\rm
Poincar\'e rank} of the system (\ref{sinv}) at $a$ is the integer $r$ of (\ref{rP}). The {\rm true Poincar\'e rank}  of (\ref{sinv}) at $a$ is  the smallest Poincar\'e rank of a local system in the meromorphic equivalence class of (\ref{sinv})  at $a$.
\end{definition}
\medskip
\noindent
We recall that the singular point $a$ is called {\it regular singular}
 if all solutions of (\ref{sinv})  have  an at most polynomial
growth as $z$ tends to $a$ in some sector with
vertex  $a$ (note that these are in general multivalued functions).
In the opposite case the singular point is called {\it irregular}.
The system (\ref{sinv}) is called {\it fuchsian} at $a$ if $r=0$, that is, if  the coefficient
form  of the system has a simple pole at $a$.

\medskip
\noindent
Assume now that $a$ is irregular. Then,  in addition
to the local monodromy matrix $G$ and the Poincar\'e rank $r$, one can 
attach {\it local Stokes data} to the system   at the singular point~$a$. These are determined as follows.
In a neighborhood of $a$ it is well-known (\cf \cite{BJL1}) that there exists a formal fundamental
solution $\hat{Y}$ of (\ref{sinv}) of the form
\begin{equation}   \label{fsol}
\hat{Y}(t)= \hat{F}(z)H(z)
\end{equation}
where $\hat{F}$ is a formal meromorphic matrix series in $z$ (in general divergent) and
$$
H(z)=(z-a)^{\tilde{J}}Ue^{Q(z)},
$$
where $Q(z)$, $U$, 
$\tilde{J}$ are block-diagonal matrices with  diagonal blocks $Q_{j}(z)$, $U_j$,  $\tilde{J}_{j}$ respectively,  $j=1,\ldots,N_{Q}$, of the same  size. We call these matrices superblocks, since they too are  block-diagonal of the  form
$$ Q_{j}(z) = \mathrm{diag}( q_j(t)I_{s_j}, q_j( t\zeta_j )I_{s_j}, \ldots,
q_j(t\zeta^{p_j -1} )I_{s_j} ), \leqno(Q)
$$
 where $q_j$ is
a polynomial in  $t=(z-a)^{-1/p_j}$ and $\zeta_j=e^{{2i\pi}/{p_j}}$, for some  integer $p_j$  not greater than
the least common multiple of  $2,3,\ldots,p$,
$$ \tilde{J}_{j} = \mathrm{diag}( J_{ s_{j} }, J_{ s_{j} } + ( 1/{ p_{j} }
)I_{s_j}, \ldots, J_{s_j} + ( (p_{j}-1)/{p_{j}} )I_{s_j} ), \leqno(J)$$
and the matrix $U_j$  decomposes into blocks $U_{lk}$ of the form
$$ U_{lk} = [ \zeta_j^{(l-1)(k-1)}I_{s_j} ] \ \ , \ \ 1 \leq l, k \leq
{ p_{j} }. \leqno(U)$$
The polynomial $q_j(t)$ has no constant term and
the integer $s_j$ is the ``multiplicity'' with which $q_j$ together with its
analytic continuations around $a$ occur on the diagonal of $Q$. As usual $I_{s_j}$
denotes the ${s_j}$-dimensional identity matrix and $J_{s_j}$ 
a constant ${s_j}$-dimensional matrix in  canonical Jordan form whose
eigenvalues $(\rho_{j}^m)_{1\leq m\le s_j}$ satisfy for all $m$ the condition $$0 \leq \mbox{Re\,}
\rho_{j}^m < 1/{p_{s_j}}. \leqno (\rho)$$

\medskip 
\noindent
In the generic case, the form of the formal fundamental matrix
has a simpler form. All superblocks $H_j$ in the decomposition of $H(z)$ are then  usual
blocks with $p_j=1$, $U=I$, that is, $H$ decomposes into a direct sum of diagonal blocks  
\begin{equation}    \label{notr}
H_j=(z-a)^{J_{s_j}}e^{q_j(z)I_{s_j}} 
\end{equation}
where $q_j(z)$ is a  polynomial of degree non greater than $r$ in $1/(z-a)$ with no constant term (with at least one $q_j$ of  degree exactly $r$) and  $J_{{s_j}}$ is a matrix in Jordan normal form
with eigenvalues  $\rho_{j}^m$ satisfying $0 \leq \mbox{Re\,}
\rho_{j}^m < 1/{p_{s_j}}$ for all $m$. 
This in particular occurs when the eigenvalues of the leading term $B_{-r-1}$ in the expansion
(\ref{rP}) are distinct. 

\noindent We will refer
to the generic case above as to the  {\it unramified} case, or case of a {\it singularity  without roots} (to the {\it ramified case} or case of {\it a singularity with roots} else). We will more precisely  say that {\it the solution} (\ref{fsol}) {\it is unramified}  if $Q$ is a polynomial in $1/(z-a)$, and that it is a {\it ramified solution} if $Q$ is polynomial in $1/t$ where $t$ is a root of $(z-a)$.  Note that $r$  in the unramified case is the  true Poincar\'e rank of (\ref{sinv}). In the general (possibly ramified) case, the true Poincar\'e rank is the least integer greater or equal to the rational degree of $Q$, that is, to the Katz rank of (\ref{sinv}) at $a$. Note that {\it via} a local meromorphic transformation it is always possible to reduce the Poincar\'e rank to the true Poincar\'e rank (for a review of general facts about the rank at an irregular singularity and rank reduction, we refer to \cite{LR} and \cite{BCM}). 

\subsection{Stokes data}

With notations as before, consider a formal fundamental solution
\begin{equation}   \label{cfsol}
\hat{Y}(z)= \hat{F}(z)(z-a)^{\tilde{J}}Ue^{Q(z)}
\end{equation}
of (\ref{sinv}) at $a$, 
where in particular $Q$ is a diagonal polynomial matrix in $1/(z-a)$ of degree $r$ with no constant term, which we  call the {\it exponential part} of $\hat{Y}$.

\medskip
\noindent
The {\it formal monodromy} (matrix) is defined as
$$
{\hat{G}}=U^{-1}\exp(2 i \pi \tilde{J})U
$$
or equivalently by
${\hat{Y}}_{2i\pi}={\hat{Y}}{\hat{G}}$
where ${\hat{Y}}_{2i\pi}$ denotes the fundamental solution obtained from $\hat{Y}$ by the change of sheet (on the Riemann surface of the
logarithm) induced by meromorphic continuation around $a$ one time in the positive (counterclockwise) direction. 

\medskip
\noindent
Note that  $Q$ and $\hat G$ are formal invariants of the system (\ref{sinv}), depending on its formal meromorphic equivalence class only.

\medskip
\noindent
Let $l_1\prec \ldots \prec l_N$ denote the {\it singular rays} of $Q$, that is, the rays from $a$ (in an affine chart containing $a$),
labeled in ascending order with respect to the positive orientation of a circle centered at $a$,  on which some  $e^{q_j - q_l}$
has maximal decay. The general theory of summability ensures that $\hat Y$ is (multi)-summable along any non-singular ray $l$
(\cf  \cite{MR}, \cite{BBRS}, \cite{Ba1}, \cite{Ba3}). If  all the polynomials $q_i-q_j$ have the same degree $r$ (this  is the case of one-level  summability)
this means that for any open sector $S$ with vertex  $a$, with opening $>\pi/r$ and bisected by   $l$,
there is a unique analytic 
fundamental matrix $Y_l$ called the {\it sum}, or $r$-{\it sum } in this case,  of $\hat Y$ along $l$, such that $YH^{-1}$ is
Gevrey $1/r$-asymptotic to $\hat F$ on this sector,  that is, for any proper subsector $S'$ of $S$ there are constants $A$ and
$C$ such that if we write

$$
\hat{F}(z)=\sum_{k=-s}^{\infty}F_k(z-a)^k,
$$
then for any $m\ge -s$

$$
| Y_j(z)H^{-1}(z)-\sum_{k=-s}^{m}F_k(z-a)^k | < CA^m(m!)^{\frac{1}{k}}(|z-a|^{m+1})
$$
as $z$ tends to $a$ in $S'$.

\medskip
\noindent
Given a  singular ray $l_i$ of $Q$, let $l_i^{-}$ and $l_i^{+}$ be two rays such that  $l_i^{-}\prec l_i\prec l_i^{+}$ and such that $l_i$ is the only singular ray contained in the oriented sector $[l_i^-,l_i^+]$. Let $Y_i^-$ and $Y_i^+$ denote the sums of $\hat {Y}$ along $l_i^-$ and $l_i^+$ respectively. Comparing these solutions on a neighbourhood of $l_i$ (they are both defined on sectors large enough to contain $l_i$) we define the {\it Stokes matrix}  with respect to the singular ray $l_i$ to be the constant matrix $C_i$ depending on $l_i$ only, such that $Y_i^-=Y_i^+C_i$. The Stokes matrices of (\ref{sinv}) at $a$ have the following properties, which we  call  {\it Stokes conditions:}
\begin{itemize}
\item For each $j$ the matrix $e^{Q(z)}C_je^{-Q(z)}$ is asymptotic to the identity
matrix $I$,
\item $C_1\cdot\ldots\cdot C_N\cdot\hat{G}=G.$ 
\end{itemize}

\noindent
Note that the first condition iin particular implies that the Stokes matrices are unipotent. The second condition is often called the {\it cyclic relation}.

\medskip
\noindent
Thus,  we have attached to each singular point $a$  of the given system,  the following data:

\medskip
\noindent
- the Poincar\'e rank at $a$ 

\smallskip
\noindent
- the proper monodromy matrix $G$ (image by (\ref{rep}) of the elementary loop $\gamma$ enclosing the singular point $a$ only)

\smallskip
\noindent
-  Stokes data which consist of the exponential part $Q$ of a formal fundamental solution (\ref{cfsol}), the formal monodromy $\hat G$ and the Stokes matrices $C_1, \ldots,C_N$ corresponding to the respective singular directions $l_1\prec\ldots \prec l_N$ of $Q$.

\medskip
\noindent
These data over all singular points of (\ref{sinv}) constitute what we will call  {\it generalized monodromy data}.

\section{The Generalized Riemann-Hilbert Problem} 
We will now give the precise terms of the generalized Riemann-Hilbert problem, GRH-problem for short.

\subsection{Statement of the problem}
We first  define the data for the inverse problem under consideration in terms of local invariants.

\begin{definition} A {\rm reduced datum} $\mathcal M$ consists of
\begin{itemize}
\item a finite subset $\mathcal{D} = \{ a_{1},\ldots,a_{n} \}$ of $\mathbb{P}^{1}(\mathbb{C})$
\item for some fixed $z_0\in{\mathbb{C}\setminus\mathcal D}$, a representation
\begin{equation}   \label{repre}
\chi:\pi_1(\mathbb{P}^1(\mathbb{C})\setminus\mathcal{D};z_0)\longrightarrow
{\mathrm {GL}}(p,\mathbb C)
\end{equation}
of the fundamental group of $\mathbb{P}^1(\mathbb{C})\setminus {\mathcal D}$ defined by matrices $G_i=\chi(\gamma_i)$ for each elementary loop-class $\gamma_i$ around $a_i, i=1,\ldots,n$,
\item for each  $a_i$ {\it local Stokes data} consisting of :

- a non-negative integer $r_i$,

- a diagonal polynomial matrix $Q_i$  in some root of $1/z_i$  (where $z_i$ denotes  a local parameter at $a_i$)  with no constant term,  with a block-diagonal decomposition $Q_i={\diag}(Q_{i,1},\dots,Q_{i,N_{Q_i}})$ in blocks of the form $(Q)$ above, and such that the fractional degree $s_i$ of $Q$ in $1/z_i$ satisfies $r_i=-[-s_i]$,

- an invertible constant matrix ${\hat G}_i$, or equivalently a matrix $\tilde{J}_i$ such that ${\hat{G}_i}=U_i^{-1}\exp(2 i \pi \tilde{J}_i)U_i$, where $U_i$ is decomposed in blocks of the form $(U)$ and $\tilde{J}_i$  is a  Jordan, block-diagonal matrix with blocks of the form $(J)$ and eigenvalues satisfying the condition $(\rho)$, and where the  size of the blocks and superblocks is determined  as before by the size of the corresponding block-decomposition of $Q_i$.

- a set of matrices $C_{i}^1,\ldots,C_{i}^{N_i}$ attached to the singular directions $l_i^1\prec\ldots\prec l_i^{N_i}$ of $Q_i$ and satisfying the Stokes conditions defined in section 2.
\end{itemize}
\end{definition}

\noindent The GRH-problem asks for the existence of a system (\ref{sinv}) with $\mathcal D$ as its set of singular points and with  $\mathcal M$ as its corresponding set of generalized monodromy data, that is, a system (\ref{sinv}) such that 

\smallskip
\noindent
 -the representation $\chi$ is the monodromy representation of (\ref{sinv}) with respect to some fundamental solution defined in the neighbourhood of $z_0$,

\smallskip
\noindent
- the system (\ref{sinv}) has Poincar\'e rank $r_i$ at $a_i$ for all $i=1,\ldots,n$,

\smallskip
\noindent
- at each $a_i$ there is a formal fundamental solution ${\hat Y}_i$ of the form (\ref{cfsol}) with $Q_i$ as its exponential part, ${\hat G}_i$ as its formal monodromy and the $C_i^j$ as its Stokes matrices along the singular lines $l_i^j$ of $Q_i$,  $j=1,\ldots,N_i$.

\bigskip
\noindent
In the following  cases, the GRH-problem reduces to classical problems.
\begin{description}
\item[(RH)]   If all  Poincar\'e ranks $r_i$ equal zero, then the data $\mathcal M$ reduce to the representation (\ref{repre}), that is,  the GRH-problem seeks a fuchsian system of linear differential equations with given singular points and a given monodromy representation. This is the classical Riemann-Hilbert problem (Hilbert's 21st problem for fuchsian equations).
\item[(BSF)] Consider the case of two singularities only, at $a_1=\infty$ and $a_2=0$, with  $\mathcal M$ data $r_1=r$, $r_2=0$ and any Stokes data at $a_1$. This is the Birkhoff standard form problem. 
\end{description}

\medskip
\noindent
{\bf Remark} The  definition of a reduced datum $\mm$ implies that  a solution to the GRH-problem for $\mm$ has a {\it minimal} Poincar\'e rank $r_i$ (equal to the true Poincar\'e rank) at each singular point $a_i$,  as it is required for the classical Riemann-Hilbert problem. 

\bigskip
\noindent
Suppose now we are given a reduced datum $\mm$. This section is devoted to the construction of a certain family $\ee$ of vector bundles with connections that realize the local data of $\mm$. Once we have achieved the construction of $\ee$, we naturally obtain the following result. 

\begin{theorem}  \label{triv}
The generalized Riemann-Hilbert problem  has a solution for $\mm$
if at least one of the vector bundles in $\ee$ is holomorphically
trivial.
\end{theorem}

\noindent In the next sections the construction of the family $\ee$ will lead to  more precise  sufficient conditions for the problem.

\subsection{Construction of $\ee$}

To solve the GRH-problem, we first apply well-known results of Malgrange and Sibuya (\cite{Mal1}, \cite{Mal2}, \cite{Si}, see also \cite{ML}]) which guarantee the existence, for each $i=1,\ldots,n$, of a local meromorphic system  
\begin{equation}   \label{si}
dy=\omega_iy
\end{equation}
of linear differential equations in a neighborhood  of $a_i$ with the given local Stokes data. 

\medskip
\noindent
The GRH-problem then can be reformulated as follows.

\medskip
\noindent
{\it Let local systems of the form (\ref{si}) be given  in neighbourhoods $O_1,\ldots, O_n$ of $a_1,\ldots, a_n$ respectively,  such that the local monodromies (with respect to suitable fundamental solutions) of these systems
generate a  representation (\ref{repre}). Does there exist a global system (\ref{sinv}) with $\{a_1,\ldots, a_n\}$ as its set of singular points and with generalized monodromy data given by those of the local systems (\ref{si}) ?}

\medskip
\noindent
A method of solution for the GRH-problem is the following. Consider a covering of $\mathbb{P}^1(\cc
)\setminus\mathcal D$ by finitely many and sufficiently small discs $U_{n+1},\ldots, U_{N}$ and
connect each $U_i$ to the base-point $z_0$ {\it via} some path
$\eta_i$ in $\mathbb{P}^1(\cc)\setminus\mathcal D$ from $z_0$ to a given endpoint in $U_i$, $i=n+1,\ldots,N$.

\medskip
\noindent
For each nonempty intersection $U_i\cap
U_j$ consider the loop $\eta_i^{-1}\circ\delta_{ij}\circ\eta_j^{-1}$, where
$\delta_{ij}$ denotes a path in $U_i\cup
U_j$ connecting 
the endpoints in $U_i
$ and $
U_j$  of $\eta_i$ and $\eta_j$  respectively. We define the constant
function
$$
g_{ij}=\chi([\eta_i\circ\delta_{ij}\circ\eta_j^{-1}]): U_i\cap U_j\longrightarrow
{\mathrm {GL}}(p,\mathbb C).
$$
It is not difficult to see that the functions $g_{ij}$ define a
gluing cocycle,   hence a vector bundle $\hat{\mathcal F}$ of rank $p$
over $\mathbb{P}^1(\cc )\setminus\mathcal D$ with these
constant transition functions.

\medskip
\noindent
For $i=n+1,\ldots,N$ consider the  system of linear differential
equations  in $U_i$
$$
dy=\omega_{i}y,\ \ \ \ \ \omega_i=0.
$$
This is a family of compatible local systems, in the following sense.
For each nonempty intersection $U_i\cap U_j$ one has
\begin{equation}    \label{gab}
\omega_i=dg_{ij}g_{ij}^{-1}+g_{ij}
\omega_j g_{ij}^{-1}
\end{equation}
since the transition functions are constant. This defines a connection $\hat{\nabla}$ on the vector bundle $\hat{\mathcal F}$,
and the forms $\omega_i$, $i=n+1, \ldots, N,$ are as usual called the {\it local forms} of the connection.

\medskip
\noindent
If we consider another coordinate description of $\hat{\mathcal F}$ by means of equivalent cocycles $g_{ij}'=\Gamma_i^{-1}g_{ij}\Gamma_i$
 where
$\Gamma_i$, $i=n+1, \ldots, N,$ denotes  a holomorphically invertible matrix  function
in $U_i$, then the corresponding local forms of the connection
$\hat{\nabla}$ are equal to $${\omega}_i'=d\Gamma_i\Gamma_i^{-1}+\Gamma_i\omega_i\Gamma_i^{-1}.$$
By construction, the connection $\hat{\nabla}$ is
holomorphic on $\mathbb{P}^1(\cc)\setminus\mathcal D$ (since all $\omega_i=0$ are holomorphic, $i=n+1, \ldots,N$) and it has the given monodromy representation (\ref{repre}).

\medskip
\noindent
We  actually can  extend $(\hat{\mathcal F},\hat{\nabla})$ to the whole Riemann sphere  by means of the local systems $\omega_i$ defined each in the neighbourhood $O_i$ of $a_i$, $i=1,\ldots,n$. This follows
from the fact that the systems (\ref{si}) have the monodromy
prescribed by the representation (\ref{repre}). Thus, one can glue
the local systems in $O_i\setminus \{a_i\}$ determined by
(\ref{si'}) and by $(\hat{\ff},\hat{\nabla})$.

\medskip
\noindent
In terms of cocycles one needs to do the following. Consider a nonempty
intersection  $O_i\cap U_{\alpha}$ and choose a fundamental solution $Y_i$
of  (\ref{si}) in this intersection.
This solution can be written as
\begin{equation}    \label{MEi}
Y_i(z)=M_i(z)(z-a_i)^{E_i},
\end{equation}
where the matrix $M_i(z)$ is holomorphically invertible  in $O_i\cap
U_{\alpha}$  and where $E_i=(1/2\pi i) \log G_i$  and the eigenvalues
$\rho_{i}^m$ of the matrix $E_i$ are normalized as follows
\begin{equation}  \label{rho}
0\leq\mbox{Re\,}\rho_{i}^m <1.
\end{equation}
Let $g_{i\alpha}(z)=M_i(z)(z-a_i)^{E_i}$. For any other
$U_{\beta}$ that has a nonempty intersection with $O_i$ consider a
path starting from a point s in $O_i\cap U_{\alpha}$ and ending
in  $O_i\cap U_{\beta}$, moving in $O_i$ around $a_i$ (less than one turn) in the
counterclockwise direction. Let  $g_{i\beta}(z)$ denote the analytic
continuation of $g_{i\alpha}$ along this path. A simple verification shows that the set $\{g_{\alpha\beta}, g_{i\alpha}\}$ defines a cocycle for the  covering $\{O_i,U_{\alpha}\}_{{1\le i\le n }\atop{n+1\le\alpha\le N}}$. Thus, one gets
a vector bundle $\ff$ on the whole Riemann sphere.

\medskip
\noindent
It follows from the preceding construction  that all the  local systems $dy=\omega_i y$, $i=1,\ldots,N$, including the systems
in the neighbourhoods $O_i$ of $a_i$, are compatible in the sense of (\ref{gab}).
Indeed, for any $i$ and $\alpha$ such that  $O_i$ and $U_{\alpha}$
have a nonempty intersection, one has
$$
dg_{i\alpha}g_{i\alpha}^{-1}+ g_{i\alpha}\omega_{\alpha}
g_{i\alpha}^{-1}=dY_iY_i^{-1}=\omega_i,
$$
that is, we get a connection $\nabla$ on the vector bundle $\ff$ with the
given local forms $\omega_i$, $i=1,\ldots,L$, and with the given monodromy
(\ref{repre}). This is the so-called {\it canonical} extension of $({\hat {\ff}},{\hat \nabla})$ in the sense of Deligne.

\medskip
\noindent
If the vector bundle $\ff$ which we have constructed  was holomorphically trivial, then
on a holomorphic trivialization of the bundle, the connection $\nabla$ would
define a global system of linear differential equations (\ref{s}) with
the given generalized monodromy data. Thus, the inverse problem would be
solved. 

\medskip
\noindent
Indeed, the triviality of the bundle $\ff$ means that for every $O_i$,
$U_{\alpha}$  and $U_{\beta}$ with $O_i\cap U_{\alpha}\ne \emptyset$, $U_{\alpha}\cap U_{\beta}\ne\emptyset$, there exist
holomorphically invertible matrix functions
$\Gamma_i$, $\Gamma_{\alpha}$ and $\Gamma_{\beta}$ (on  $O_i$,
$U_{\alpha}$, $U_{\beta}$ respectively)  such that
$$
\Gamma_{\alpha}=\Gamma_ig_{i\alpha},\ \ \ \Gamma_{\beta}=\Gamma_{\alpha}
g_{\alpha\beta}
$$
if the corresponding intersections are not empty. This implies that
the forms
$$
\omega_i'=d\Gamma_i\Gamma_i^{-1}+\Gamma_i\omega_i\Gamma_i^{-1},\ \ \
\omega_{\alpha}'=d\Gamma_{\alpha}\Gamma_{\alpha}^{-1}+\Gamma_{\alpha}
\omega_{\alpha}\Gamma_{\alpha}^{-1}
$$
coincide over the intersections of the corresponding pieces of
the covering and thus define a global form $\omega$. The
fundamental matrices  of the new and original local systems
 are connected by gauge transformations $Y_i'=\Gamma_iY_i$, which implies that the constructed
system  has the required generalized monodromy data.

\medskip
\noindent
Unfortunately, the bundle $\ff$ as a rule is not holomorphically trivial.
But it turns out that this bundle is always meromorphically trivial.
More precisely, for any choice of a point $b$ in some $O_l$ there exists a meromorphic
trivialization of the bundle which is holomorphic outside of $\{b\}$. In terms
of a cocycle involving $O_l$,  as above, this means that the desired functions $\Gamma_i$ will be
holomorphically invertible for $i\ne l$, and 
$\Gamma_l$  meromorphic only at $b$ (and holomorphically invertible
in $O_l\setminus \{b\}$).

\medskip
\noindent
Choose $b=a_l$ for some  $l$, $1\le l\le n$.  From the meromorphic
trivialization $\{\Gamma_i\}$ of the bundle $\ff$ we get a global
system (\ref{sinv}) with all the given generalized monodromy data
except one, namely the Poincar\'e rank at $a_l$ which may be greater than the given integer $r_l$, since the matrix $\Gamma_l$ is
meromorphic only at $a_l$.  And
for a number of inverse problems such as the classical
Riemann-Hilbert problem or the problem of the standard Birkhoff
form, we must realize this datum as well.

\medskip
\noindent
To achieve this, we shall replace the local systems (\ref{si}) in the construction of $(\ff,\nabla)$ by new systems
\begin{equation}  \label{si'}
dy=\omega_i ' y
\end{equation}
with
$$\omega'_i=d\Gamma_i\Gamma_i^{-1} + \Gamma_i\omega_i\Gamma_i^{-1},$$
{\it via}  (families of) gauge transformations $y'_i=\Gamma_i y$, where $\Gamma_i$ is holomorphically invertible in $O_i\setminus\{a_i\}$ and meromorphic at $a_i$. 

\medskip
\noindent
\begin{definition}\ \ Assume  $a_i$ is a singular point at which the formal solution $\hat{Y}_i$
of (\ref{si}) is unramified. This means that $\hat{Y}_i$ has the form
(\ref{fsol}), (\ref{notr}). An {\rm
admissible matrix} is an
integer-valued diagonal
matrix
$\Lambda_i=\mbox{diag\,} 
 (\Lambda_i^1,\ldots,\Lambda_i^{N_{Q}})$, that is, a diagonal matrix whose entries are integers,  blocked in the same way  as $Q(z)$ and such that the matrix function
$$
(z-a_i)^{\Lambda_i}\tilde{J}(z-a_i)^{-\Lambda_i}
$$
is holomorphic at $a_i$. 
\end{definition}

\medskip
\noindent
Note that any diagonal integer-valued matrix $\Lambda_i$ whose diagonal elements form a non-increasing sequence  is
admissible and that the set of admissible matrices is infinite.

\medskip
\noindent
The matrix $\hat{Y}_i$ can be written as
follows
$$
\hat{Y}_i(z)=\hat{F}(z)(z-a_i)^{-\Lambda_i}(z-a_i)^{\Lambda_i}H(z)
$$
with $H(z)$ as in (\ref{fsol}). (For simplicity of notation, we will  omit the index $i$ when introducing new functions, although all  calculations  depend on $a_i$).
The formal matrix function $\hat{F}'(z)=\hat{F}(z)(z-a_i)^{-\Lambda_i}$
is still meromorphic at $a_i$. The proof of the following technical lemma
proceeds as for Sauvage's
 lemma in  \cite{Bol1}.
\begin{lemma}   \label{Sauv}
For any formal meromorphic matrix $F( \tau ) \in
\mathrm{GL}(p, \mathbb{C}[[ \tau ]] [1/ \tau])$  there exists a
matrix $\Gamma( \tau )$, polynomial in $1/\tau$ and holomorphically invertible
outside of $\tau=0$, such that
$$ \Gamma( \tau ) F( \tau ) = {\tau}^{K} F_{0}( \tau ), $$
where $K$ is a diagonal integer-valued matrix and
 $F_{0}( \tau )$ is an invertible formal holomorphic matrix series in $\tau$.
\end{lemma}
If we apply this lemma to the matrix $\hat{F}'(z)$ we get
$$
\Gamma(z-a_i)\hat{F}'(z)=(z-a_i)^K\hat{F}_0(z).
$$

\medskip
\noindent
Let us transform the local system (\ref{si}) into (\ref{si'}) via the meromorphic gauge transformation $y'=\Gamma_i y$ where 
$\Gamma_i(z)=(z-a_i)^{-K}\Gamma(z-a_i)y$.
We get a formal fundamental  solution of the new system,  of the form
\begin{equation}   \label{Ffsol}
\hat{Y}_i'(z)=\hat{F}_0(z)(z-a_i)^{\Lambda_i}(z-a_i)^{\tilde{J}}e^{Q(z)}.
\end{equation}

\medskip
\noindent
This transformation does not increase the Poincar\'e rank $r_i$.
Indeed,  the form $\omega_i'$ of  (\ref{si'}) can be written as
$$
\omega_i'=d\hat{Y}_i'(\hat{Y}_i')^{-1}=d\hat{F}_0(\hat{F}_0)^{-1}+
\hat{F}_0\frac{1}{z-a_i}\left(\Lambda_i+\right.$$
\begin{equation}   \label{calc}
+\left.(z-a_i)^{\Lambda_i}\tilde{J}(z-a_i)^{-\Lambda_i}+(z-a_i)\frac{dQ(z)}{dz}
\right)
(\hat{F}_0)^{-1}dz.
\end{equation}
The fact that the matrix $\Lambda_i$ is admissible and
the matrix $\hat{F}_0$ invertible, and  the fact that the degree of
$Q(z)$ is equal to $r_i$ (with respect to $1/(z-a_i)$) together guarantee
that the Poincar\'e rank at $a_i$ of the new local system remains
equal to $r_i$.

\medskip
\noindent
Let us replace the initial local system (\ref{si}) 
in $O_i$ with the system (\ref{calc}), which we will write
\begin{equation}   \label{sLi}
dy=\omega^{\Lambda_i}y
\end{equation}
to keep track of the
admissible matrix $\Lambda_i$ used in the construction. Let us extend the initial vector bundle
$(\hat{\ff},\hat{\nabla})$, constructed from the  representation (\ref{repre}), 
over the point $a_i$ using this new system
(instead of the initial one).

\medskip
\noindent
Assume  that $a_i$ is a regular singular point. Consider in this case
 an analytic fundamental solution $Y_i(z)=M_i(z) (z-a_i)^{E_i}$ such that
moreover the matrix $E_i$  has an upper triangular  form and the
entries $e_{kl}$ of $E_i$ equal zero if $\rho_i^k\neq\rho_i^l$,
where the complex numbers $\rho_i^m$ denote the eigenvalues of $E_i$. Since
 $a_i$ is regular singular, the matrix $M_i(z)$ 
is meromorphic at $a_i$. Thus, we can follow the same procedure as in the
case  of an irregular point without roots, to construct a new system via
an admissible matrix $\Lambda_i$, where admissibility here means that  the
matrix $(z-a_i)^{\Lambda_i}E_i(z-a_i)^{-\Lambda_i}$ is holomorphic at $a_i$.

\medskip
\noindent
Now assume  that $a_i$ is an irregular singular point  with roots.
By an admissible matrix $\Lambda_i$ we mean here a diagonal
integer-valued matrix
$\Lambda_i=\mbox{diag\,}(\Lambda_i^1,\ldots,\Lambda_i^{N_{Q}})$
blocked in the same way  as $Q(z)$ and such that the matrix function
\begin{equation}   \label{adroots}
(z-a_i)^{\Lambda_i^j}\tilde{J}_j(z-a_i)^{-\Lambda_i^j}
\end{equation}
is holomorphic at $a_i$ if the superblock $Q_j$ has no ramification, and
$\Lambda_i^j$ is a scalar matrix if the superblock $Q_j$ has ramification.

\medskip
\noindent
Let us proceed  with the system (\ref{si}) at an irregular singular point
with roots ({\it i.e.} the formal solution is ramified) in the same way as in the unramified case. Again, we get a system  (\ref{sLi}) with the same local Stokes data as the initial  one.

\medskip
\noindent
Choose a collection $\Lambda=(\Lambda_1,\ldots,\Lambda_n)$ of admissible matrices (in the above sense, depending on the type of the singularity $a_i$) and consider the extension $(\ff^{\Lambda},\nabla^{\Lambda})$
of  $(\hat{\ff},\hat{\nabla})$ over the singular points $a_i$ {\it via} the systems (\ref{sLi}) obtained by means of the matrices $\Lambda_i$. Then, by construction,  the extended connection
$\nabla^{\Lambda}$ has the given Poincar\'e ranks and generalized monodromy
data. We get in this way an infinite set $\ee$ of vector bundles
$(\ff^{\Lambda},\nabla^{\Lambda})$  with connections that have  the prescribed
generalized monodromy data. 

\medskip
\noindent It follows immediately from previous considerations that Theorem \ref{triv} holds for this family $\ee$ of vector bundles.

\medskip
\noindent
Note  that the converse of Theorem \ref{triv} is not true, since $\ee$ does not contain all   vector bundles with connections
 having the prescribed generalized monodromy data. The reason for this is
that there are local systems (\ref{si}) with the given data  whose
formal fundamental matrix cannot be written in the form
(\ref{Ffsol}) with an invertible matrix $\hat{F}_0$ (
this  in particular occurs at any regular, but not fuchsian singularity). This situation differs significantly  from the fuchsian case (where all
Poincar\'e ranks equal zero) in which Theorem \ref{triv} gives
necessary and sufficient conditions for the positive solvability
of the Riemann-Hilbert problem (see \cite{AB}, \cite{Bol1}).

\section{Sufficient Conditions for the Generalized Riemann-Hilbert Problem  }

We keep  notation from section 2. Consider a bundle $(\ff^{\Lambda},\nabla^{\Lambda})$ in  $\ee$.
It follows from (\ref{fsol}), (\ref{MEi}) and (\ref{calc}) that
$$
\mbox{tr\,}\omega^{\Lambda_i}=\mbox{tr\,}(\Lambda_i+T_i)\frac{dz}{z-a_i}
+ \ \ \ \mbox{a holomorphic form},
$$
where $T_i=\tilde{J}$ in the irregular case and  $T_i=E_i$ if 
$a_i$ is a regular singular point.
The eigenvalues $\beta_{i,m}$ of the matrix $\Lambda_i+\tilde{J}$ in the unramified scase 
(resp. of the matrix $\Lambda_i+E_i$ in the regular singular case)
are called {\it formal exponents} ({resp. \it exponents}) of the connection
$\nabla^{\Lambda_i}$  at $a_i$. 
The {\it degree} $\ {{\mathrm {deg} }\,}\ff^{\Lambda}$ of the bundle $\ff^{\Lambda}$
is by definition the sum
$$
{{\mathrm {deg}}\,}\ff^{\Lambda}=\sum_{i=1}^n
\mbox{res}_{a_i}\mbox{tr\,}\omega^{\Lambda_i}=
 \sum_{i=1}^n\mbox{tr\,}(\Lambda_i+T_i).
$$
Let us recall that a bundle $\ff$ is called  {\it stable} (respectively {\it
semistable}) if for any proper subbundle $\ff'$ of $\ff$,  the
{\it slope} $\mu(\ff')=\dd (\ff')/\rr (\ff')$ of $\ff'$ is less (resp. non greater)
than the slope  $\mu(\ff)$ of $\ff$.

\medskip
\noindent
A holomorphic bundle on the Riemann
sphere is trivial if and only if it is a semistable bundle of
degree $0$. Indeed, each   vector bundle $\ff$ on the Riemann sphere is
holomorphically equivalent to a sum of line bundles
\begin{equation}  \label{split}
\ff\cong {\cal O}(c_1)\oplus\ldots\oplus {\cal O}(c_p),
\end{equation}
where the ordered set of integers  $c_1\geq\ldots\geq c_p$ is called the {\em 
splitting type} of the bundle $\ff$. If the bundle $\ff$ is semistable  of
degree zero, then  $c_1+\ldots+c_p=0$ and $c_i\leq 0$, $i=1,\ldots,p$.
Thus, $c_1=\ldots=c_p=0$ and $\ff$ is holomorphically trivial.

\medskip
\noindent
In what follows  we will need  the notion of stability of a pair
consisting of a vector bundle and a connection. A
subbundle $\ff'$ of the bundle $\ff^{\Lambda}$ is said to be {\it  stabilized by
the connection $\nabla$} if the covariant derivative $\nabla_{{d}/{dz}}$
maps local holomorphic sections of  $\ff'$  into
sections of the same subbundle. In the coordinate  description $\{O_i, U_{\alpha}\}$, $\{g_{i\alpha}, g_{\alpha \beta}\}$, $\{\omega^{\Lambda_i}, \omega_\alpha\}$ of the pair $(\ff^{\Lambda},\nabla^{\Lambda})$, the existence of such a subbundle means the following. For all $i=1,\ldots,n$, there exist matrices $\Gamma_i$, each holomorphically invertible in the corresponding  $O_i$ and such that all systems  (\ref{si'}) obtained from the systems (\ref{sLi}) via the gauge transformations $\Gamma_i$ have the form
\begin{equation}    \label{redomega}
\omega_i'=\left(\begin{array}{cc}
\omega_i^1 & *\\
0 & \omega_i^2
\end{array}\right)
\end{equation}
with  blocks $\omega_i^1$ of the same size for all $i$.
The  local subsystems $\omega_i^1$ define the
restriction of the connection $\nabla^{\Lambda}$ to a subbundle
$\ff'$. Each formal solution $\hat{Y}_i'$ of such a
system (\ref{si'})  can be chosen to have
the same upper block-triangular structure
$$
{\hat{Y}}_i'=\left(\begin{array}{cc}
{\hat{Y}}_i^1 & *\\
0 & {\hat{Y}}_i^2
\end{array}\right),
$$
where the matrix ${\hat{Y}}_i^1$ serves as a formal fundamental matrix for
the subsystem  $\omega_i^1$. Moreover, the matrix $\hat{Y}_i'$ is
connected to the initial matrix ${\hat{Y}}_i$ by
$\hat{Y}_i'=\Gamma_i{\hat{Y}}_iS$ where $S$ is a constant invertible matrix and $\hat{Y}_i'$  has the same form (\ref{Ffsol})
as the  initial matrix ${\hat Y}_i$, namely
\begin{equation}   \label{redfsol}
\hat{Y}_i'(z)=\hat{F}_0'(z)(z-a_i)^{\Lambda_i'}(z-a_i)^{\tilde{J}'}e^{Q'(z)},
\end{equation}
where   $\Lambda_i'=S^{-1}\Lambda_iS$, $\tilde{J}'=S^{-1}\tilde{J}S$,
$Q'(z)=S^{-1}Q(z)S$, where the  matrices  $\Lambda_i'$ and  $Q'(z)$
are  diagonal and obtained by suitable permutations of the diagonal elements
of $\Lambda_i$ and  $Q(z)$ respectively, and where the matrix $\tilde{J}'$ is upper triangular and
satisfies the following condition: the entry $e_{kl}$ of the matrix equals
zero if for corresponding eigenvalues of $\tilde{J}'$ one has
$\rho_i^k\neq\rho_i^l$. Moreover the invertible formal holomorphic matrix
$\hat{F}_0'(z)$ has the same upper block triangular structure as the matrix
$Y_i'$. The existence of such a fundamental matrix follows from results
of \cite{BJL3} and \cite{AB}.

\medskip
\noindent
Thus, the degree of the subbundle $\ff'$ can be
determined   in the same way as the degree  of $\ff^{\Lambda}$, after replacing 
the systems (\ref{sLi}) by the subsystems $\omega_i^1$ of (\ref{redomega}).

\medskip
\noindent
Let us recall that a pair $(\ff,\nabla)$ consisting of a vector bundle $\ff$ and a connection
$\nabla$ on $\ff$ is said to be {\it stable}  (resp.{\it semistable}) if for any proper subbundle
$\ff'$ of  $\ff$ stabilized by $\nabla$  one has
$\mu(\ff')<\mu(\ff)$ (resp. $\mu(\ff')\leq\mu(\ff)$) (\cf  \cite{Sim}).

\medskip
\noindent
The main result of our paper is the following.
\begin{theorem}   \label{main}

Let $\mm$ be a reduced datum, in which one at least of the prescribed singularities is without roots.
If there exists a collection
$\Lambda$ of admissible  matrices $\Lambda_1,\ldots,\Lambda_n$ such that the
corresponding pair $(\ff^{\Lambda},\nabla^{\Lambda})$ is stable,
then the generalized Riemann-Hilbert problem for $\mm$  has a solution.\end{theorem}
\noindent
Theorem \ref{main} is analogous to Theorem 1 of \cite{Bol5} which was proved by A. Bolibrukh in the case of the classical Riemann-Hilbert problem. The proof below follows the proof of \cite{Bol5} with some simplifications (Theorem 1 of  \cite{Bol5} was proved for any compact Riemann surface).

\begin{proof}
Without loss of generality we may assume that $i=1$ and $a_1=0$.
We will denote a stable pair $(\ff^{\Lambda},\nabla^{\Lambda})$ in
$\ee$ simply by $(\ff,\nabla)$.

\smallskip
\noindent
We first consider the case in which all  eigenvalues of the
local monodromy matrices (formal and proper) are positive
real numbers. The real parts of all   $\rho_{i,m}$ then equal zero. The
pair $(\ff,\nabla)$ being stable,  one has
$\mu(\ff ')<\mu(\ff)$ for every (proper, nonzero) subbundle $\ff '$ stabilized by the connection $\nabla$ . All   $\mu(\ff ')$
are rational numbers whose denominators are not greater than
the rank of $\ff$,  hence the set of such numbers is finite
and the number  $\mu_{max}=\max_{\ff '}(\mu(\ff ')-\mu(\ff))$ is
well defined. The stability of the pair $(\ff,\nabla)$ implies
 $\mu_{max}<0$.

\smallskip
\noindent
 Starting with the initial pair $(\ff,\nabla)$, we construct
a new pair pair $(\bar{\ff},\bar{\nabla})$ as follows.

\noindent
We replace the matrices $\Lambda_i$ which were used to construct the initial
bundle $\ff$, with  matrices $\Lambda_i'=N\Lambda_i$
for some   positive integer  $N$. We choose $N$ such that $-N\mu_{max}>>(R-2+n)p^3$,
where $R$ is the sum of the Poincar\'e ranks of $\nabla$ over all singular points, and where $>>$  means  ``sufficiently larger than" (the difference can be made as large as needed). The new matrices $\Lambda_i'$
are clearly admissible and thus, the corresponding vector bundle  $\bar{\ff}$
is equipped with a connection $\bar{\nabla}$ which has the same generalized
monodromy data (including the Poincar\'e ranks at the
singularities) as the initial connection.

\medskip
\noindent
It is not difficult to see that
the pair $(\bar{\ff},\bar{\nabla})$ is stable, and that for
any subbundle
$\bar{\ff}'\subset \bar{\ff}$ stabilized by the connection
$\bar{\nabla}$ one has
\begin{equation}   \label{<<}
\mu(\bar{\ff}')-\mu(\bar{\ff})\leq N\mu_{max}<<-(R-2+n)p^3.
\end{equation}
Indeed, if $\ff '\subset \ff$ is a given (proper, nonzero) subbundle of $\ff$ that is stabilized by $\nabla$, then, as explained earlier, each local form
$\omega_i'$ of the connection $\nabla$ (in the corresponding
coordinate description of the bundle $\ff$) is of the form (\ref{redomega})
at  $a_i$, where the first diagonal block $\omega^i_1$
corresponds to $\ff '$ and  the sum of traces of
$\mbox{res}_{a_i}\omega^i_1$ over all singular points $a_i$ is equal to the degree
of $\ff '$. It follows from (\ref{redfsol}) that replacing $\Lambda_i$ by $\Lambda_i'$   multiplies
the real parts of the traces by $N$. Thus, all  degrees of
all proper subbundles stabilized by $\nabla$  are multiplied by $N$.
(Note that  the sum of the imaginary parts of  the corresponding
traces for every stabilized subbundle is zero; since we assumed
that the eigenvalues of all  local monodromy matrices are positive
real,  only  admissible matrices $\Lambda_i$ will give a real input in
the degrees.)  The slopes of all stabilized subbundles (including the bundle $\ff$) are  multiplied by $N$, hence 
$$\mu(\bar{\ff }')=N\mu(\ff ')<<\mu(\bar{\ff})=N\mu(\ff),$$  
and the new pair is stable.

\medskip
\noindent
If the entries $\lambda_{i,k}$ and $\lambda_{i,l}$ of the matrix
$\Lambda_i$ are distinct, then after multiplication by~$N$
 their difference will be
``sufficiently larger" than $(R-2+n)p^3$.

\medskip
\noindent
Since  further in the  proof we shall need  ``large'' differences between
 entries of the  matrix $\Lambda_1$  at $a_1=0$, we
have to modify  the matrix $\Lambda_1'$ once more to separate
possible pairs of equal eigenvalues. In order to preserve the stability of the pair,  replace the matrix
$\Lambda'_1=\mbox{diag}\,(\lambda'_1,\ldots,\lambda'_p)$ with an
admissible matrix $\Lambda^{''}_1=
\mbox{diag}\,(\nu_1,\ldots,\nu_p)$ such that
$\mbox{tr}\,\Lambda^{''}_1=\mbox{tr}\,\Lambda'_1$ and such that
$$ |\nu_i-\lambda_i'|< (R-2+n)p^3,$$
\begin{equation}   \label{''}
(R-2+n)p^2\leq \nu_i-\nu_{j}<(R-2+n)p^3 \ \ \mbox{for}\  i<j 
\end{equation}
if $\lambda_i'=\lambda_j'$.
To do this it is sufficient to replace
every maximal chain of equal numbers
$\lambda_{i_1}'=\ldots=\lambda_{i_s}'$, $i_1<\ldots<i_s$ with the chain
$$
\nu_{i_t}=\lambda_{i_t}'+(R-2+n)p^2 ([s/2]-t+1) \mbox{    for  } t\leq [s/2],
$$
$$
\nu_{i_t}=\lambda_{i_t}'+(R-2+n)p^2 ([(s+1)/2]-t) \mbox{    for  } t> [s/2],
$$
where $[\ ]$ stands for the integer part.
Let  $\Lambda_1$  denote again this admissible matrix $\Lambda_1''$. Let $\gg$ denote the corresponding bundle
and $\nabla$ the corresponding logarithmic
connection on $\gg$. From the construction (in particular from (\ref{<<}))
it follows that the pair $(\gg,\nabla)$ is stable.

\medskip
\noindent
Assume that the bundle $\gg$  is non trivial and  consider a meromorphic
trivialization  of $\gg$, holomorphic outside of $a_1=0$.
As  was explained in section 2, the corresponding global system
(\ref{sinv}) constructed {\it via} this trivialization has the prescribed generalized
monodromy data except at $a_1=0$, where the Poincar\'e rank may be greater than
$r_1$. This means that the formal fundamental matrix $\hat{Y}_1(z)$ of the
system
is of the form (\ref{Ffsol})
$$
\hat{Y}_1(z)=\hat{F}(z)z^{\Lambda_1}z^{\tilde{J}}e^{Q(z)},
$$
where the  matrix $\hat{F}$ is a formal meromorphic series (and not
formal holomorphic as it would  be if $\gg$ were holomorphically
trivial). Moreover, in view of the decomposition (\ref{split}) one can choose a
meromorphic trivialization of the bundle (that is,  matrices $\Gamma_i$)
such that the matrix $\Gamma_1(z)$ is of the form
$$
\Gamma_1(z)=z^{-K}\Gamma_1^0(z),
$$
where $K$ is the integer-valued diagonal matrix $K={\mathrm {diag}}(c_1,\ldots,c_p)$, $c_1\ge\ldots \ge c_p$,  and $\Gamma_1^0(z)$ is holomorphically invertible in
$O_1$. Thus, the formal fundamental matrix $\hat{Y}_1$ is of the form
\begin{equation}   \label{rfsol}
\hat{Y}_1(z)=z^{-K}\hat{F}^0(z)z^{\Lambda_1}z^{\tilde{J}}e^{Q(z)}
\end{equation}
with an invertible formal holomorphic series $\hat{F}^0(z)$.

\medskip
\noindent
The following statement generalizes related
results of \cite{Bol5} and  plays a crucial role in
the proof.
\begin{lemma}   \label{estim}
The following inequalities hold
for the entries of the matrix $K$ in the decomposition (\ref{rfsol}):
$$
c_j-c_{j+1}\leq R+n-2,\ \ \ \ \  j=1,\ldots,p-1,
$$
where  $R$ denotes the sum of the (prescribed) Poincar\'e ranks $r_i$
of the
initial connection $\nabla$ at all singularities.
\end{lemma}
\begin{proof}
Assume that  $c_l-c_{l+1}>R+n-2$ for some $l$. This contradicts the fact that the pair
$(F,\nabla)$ is stable, hence prove the lemma.

\bigskip
\noindent
In view of (\ref{rfsol}) the form $\omega$, for the system (\ref{sinv}) constructed above, can be written 
as follows in $O_1$,
$$
\omega=z^{-K}\theta z^K,
$$ where the form $\theta$ has a pole of order $r_1+1$ since
$$
\theta= -{\frac{K}{z}} + d\left(\hat{F}^0(z)z^{\Lambda_1}z^{\tilde{J}}e^{Q(z)}\right)
\left(\hat{F}^0(z)z^{\Lambda_1}z^{\tilde{J}}e^{Q(z)}\right)^{-1}
$$
following the calculation in (\ref{calc}). 

\medskip
\noindent
The entries $\omega_{mj}$ and
$\theta_{mj}$ of the matrix
differential forms $\omega$ and $\theta$  respectively, for $m\neq j$,  are
connected as follows
$$
\omega_{mj}(z)=\theta_{mj}(z)z^{-c_m+c_j},
$$
and we have by assumption $c_j-c_m>R+n-2$ for $m>l,\ j\leq l$.
The orders of zero of the differential forms $\omega_{mj}(z)$  at
$a_1=0$, for $m>l,\ j\leq l$,  are therefore greater than  $R+n-r_1-3$, whereas the sum of the orders of
poles at the other singular points is not greater than
$R-r_1+n-1$ (respectively $R-r_1+n-3$) if the point at infinity is non-singular (resp. singular). If the form
$\omega$ is holomorphic at infinity, then it has a zero of order
two there. One gets in both cases that for each entry
$\omega_{mj}(z)$, $m>l,\ j\leq l$,  the degree of its singular
divisor (the sum of orders of zeros minus the sum of orders of poles
on the Riemann sphere) is greater than zero.
Thus, all such entries $\omega_{mj}$ equal zero identically  and  
$\omega$ has the form
\begin{equation}   \label{red}
\omega'=\left(\begin{array}{cc}
\omega^1&\ast\\
0 & \omega^2
\end{array}\right),
\end{equation}
where the  form $\omega^1$ has size $l\times l$.

\medskip
\noindent
This implies that there exists a constant invertible matrix $S$ such
that $\hat{Y}_1(z)S$ has a form similar to
(\ref{rfsol})
$$
\hat{Y}_0(z)=\hat{Y}_1(z)S=z^{-K}\hat{F}^0(z)z^{\Lambda_1'}z^{\tilde{J}'}
e^{Q'(z)},
$$
where
$$
\hat{F}^0= \left(\begin{array}{cc}
\hat{F}^1&\ast\\
0 & \hat{F}^2
\end{array}\right),
$$
and $\hat{F}^1$ is of size $l\times l$.
The vector bundle  $\ff^1$ of rank $l$ carrying the connection
$\nabla^1$ defined by the subsystem $\omega^1$ is a subbundle of
$\gg$ which is stabilized by  $\nabla$. The degree of
this subbundle is $c_1+\dots+c_l$. It follows from the assumption $c_l>c_{l+1}$ that
$$
\mu(\ff^1)=\frac{c_1+\ldots+c_l}{l}>\frac{c_1+\ldots+c_p}{p}=
\mu(\gg).
$$
This contradicts the semistability of $(\gg,\nabla)$, hence it proves the lemma.

\end{proof}

\noindent
Notice that we have so far only  used the semistability  of the pair
$(\gg,\nabla)$ (which is weaker than its stability). In terms of vector
bundles the previous lemma can be reformulated as follows.

\begin{lemma}   \label{semi}
If a pair $(\gg,\nabla)$ is semistable, then  the inequalities 
$$
c_j-c_{j+1}\leq R+n-2,\ \ \ \ \  j=1,\ldots,p-1
$$ hold for the splitting type $c_1\ge\ldots\ge c_p$
of  $\gg$ and for the sum
$R$  of  all  Poincar\'e ranks of the
connection $\nabla$ at the $n$ prescribed singular
points.
\end{lemma}
\medskip
\noindent
Let us return to the proof of the theorem.  We will also need the
following technical  lemma, which is given in  \cite{Bol1} and \cite{Bol3}.
\begin{lemma}\label{TL}
Let the matrix $\hat{F}(z)$ (formal or analytic) be  invertible in a
neighborhood
$O_1$ of $a_1=0$. Then for any integer-valued diagonal
matrix $K={\mathrm {diag}}(k_1,\ldots , k_p)$ there exists a matrix $T(z)$,
polynomial in $1/z$ and holomorphically invertible outside of  $\{a_1\}$, such that
\begin{equation} \label{Com}
T(z)z^{K}\hat{F}(z)=\hat{H}(z)z^{D},
\end{equation}
where the matrix $\hat{H}(x)$ is   invertible in $O_1$
and $D$ is a diagonal matrix obtained by a suitable permutation of
the diagonal elements of $K$.
\end{lemma}

\medskip
\noindent
Apply this lemma to the factor $z^{-K}\hat{F}^0(z)$ in the expression (\ref{rfsol}) of the
fundamental matrix $\hat{Y}_1(z)$. The gauge
transformation $\hat{Y}_1^f(z)= T(z)\hat{Y}_1(z)$ (which is
holomorphically invertible outside of zero) changes our system
(\ref{sinv}) into a system with the following formal fundamental
matrix at $a_1=0$
$$
\hat{Y}_1^f(z)=\hat{H}(z)z^{D+\Lambda_1}z^{\tilde{J}}e^{Q(z)},
$$
where the formal matrix series  $\hat{H}(z)$ is invertible. It follows from Lemmas
\ref{estim} and \ref{TL}  that the difference between  any
two diagonal elements of  $D$ is bounded by 
$(R+n-2)(p-1)$. Since by construction the matrix
$\Lambda_1$ satisfies the inequalities (\ref{''}),  the
diagonal entries of  $D+\Lambda_1$ form a decreasing
sequence; hence this matrix is admissible. From
(\ref{calc}) we deduce that the final system has Poincar\'e rank $r_1$ at
 $a_1=0$. And since $T$ is holomorphically
invertible outside of $\{a_1\}$ we get that the final system has the
required Poincar\'e ranks at all points. The theorem is proved
(under the assumptions made at the beginning of the proof).

\medskip
\noindent
Now consider the  case of arbitrary eigenvalues of the local monodromy
operators. For any $N\in \mathbb Z$, we have
$$
N\beta_i^j=N\lambda_i^j+N\rho_i^j=N\lambda_i^j+(N-1)\rho_i^j +\rho_i^j=
{\lambda '}_i^j +\rho_i^j +\alpha_i^j,
$$
where ${\lambda'}_i^j=N\lambda_i^j+[\mbox{Re\,} (N-1)\rho_i^j]$
and $0\leq\mbox{Re\,}\alpha_i^j<1$  (where   $[\ ]$ as before stands
for the integer part).

\medskip
\noindent
For each $i=1,\ldots,n$, consider formal solutions
${\hat Y}_i={\hat F}_iH_i$, where ${\hat F}_i = z^{\Lambda_i} z^{\tilde J} e^Q$, and replace the matrix $\Lambda_i$ with ${\Lambda'}_i$ obtained as follows. Let  $\Lambda_i^j$ denote the blocks of $\Lambda_i$ corresponding to the superblocks $H_i^j$ of $H_i$.  If $H_i^j$ has no ramification, replace the block $\Lambda_i^j$ with the diagonal matrix $ '\Lambda_i^j$ with entries ${\lambda'}_i^j$. If $H_i^j$ is ramified, of size $d$, then replace the block $\Lambda_i^j$ with ${}^N\Lambda_i^j=N\Lambda_i^j+sI_d$, where
$$
s=\left[\frac{\sum \mbox{Re\,} (N-1)\rho_i^j}{d}\right],
$$
and where the sum is taken over all eigenvalues of the
superblock $\tilde{J}_j$.  The matrices $\Lambda_i'$ are clearly admissible. 

\medskip
\noindent
Let us prove  that for sufficiently large $N$ the corresponding pair
$(\gg',\nabla')$ is stable and the inequalities (\ref{<<}) hold.
If for each $i$ the formal solution at $a_i$ is unramified, then
this is clear from the construction, since the sum of all $\mbox{Re\,}\alpha_i^m$, 
$0\leq\mbox{Re\,}\alpha_i^m<1$, over all $i,m$, is bounded by
$pn-1$. Thus, in the first change of degrees carried out at the
beginning of the proof one only needs to   replace all the  slopes
$\mu_l$ (obtained after multiplication by $N$) of all subbundles
stabilized by the connection,  by numbers $\mu_l+t_l$, where $|t_l|$
is bounded by $pn-1$, which is $<< N\mu_{max}$.

\medskip
\noindent
Let for some $i$ the formal solution have a superblock $H_j(z)$
with ramification. From the related result of \cite{BJL1}  it
follows  that if our bundle has a subbundle stabilized by the
connection and if the local form of the connection has the form
(\ref{redomega}), then the part of the formal solution at $a_i$ that corresponds to the
superblock $H_j$ appears entirely as a block of the
 formal solution of either   the subsystem $\omega_i^1$ or
 the corresponding quotient system
$\omega_i^2$ (in \cite{BJL1} this property is called  {\it
irreducibility of superblocks with roots}).

\medskip
\noindent
Thus, when we replace a block $\Lambda_i^j$ with ${}^N\Lambda_i^j$
instead of ${}'\Lambda_i^j$,  we actually replace all the final slopes  $\mu_l$
of  subbundles stabilized by the connection, with numbers
$\mu_l+t'_l$  where $|t'_l|$ is not greater than $p(n+1)-1$, which again
is $<< N\mu_{max}$.

\medskip
\noindent
Thus,  the  pair which we have constructed iis stable and the inequalities (\ref{<<})
hold. The remainder of the proof for  the general case
is  the same as in the special case considered before (of monodromy operators with positive eigenvalues).
\end{proof}

\begin{definition}  \label{gendat} If $\mm$ is a reduced datum, let $\mm_s$  denote the  reduced datum consisting  of $\mm$ and of a family of local systems (\ref{si}) realizing  $\mm$. A datum $\mm_s$ is said to be $\mathrm{generic}$ if the pair $(\ff,\nn)$, or canonical extension, constructed as before from the systems (\ref{si}),  has no subbundle stabilized by $\nn$.
\end{definition}

\noindent
 In terms of local systems,  genericity  means that it is impossible to transform the systems (\ref{si}) in the form (\ref{redomega}) by means of local
holomorphic gauge transformations.

\medskip
\noindent
Any pair
$(\ff^{\Lambda},\nabla^{\Lambda})$ constructed from a generic datum $\mm_s$ is
by definition stable,  and we obtain the following statement which generalizes a
similar result of \cite{Bol2}.

\begin{corollary}   \label{generic}
Let $\mathcal M_s$ be a generic  datum.  Then the
generalized Riemann-Hilbert problem for $\mathcal M$ has a 
solution if one at least of the prescribed singularities is without roots.
\end{corollary}

\noindent
If in particular the monodromy representation (\ref{repre}) is irreducible, then
$\mathcal M$ is clearly generic;  thus we obtain the expected generalization of 
\cite{Bol1} and \cite{Ko}.

\begin{corollary}   \label{irr}
Let  $\mm$ be a reduced datum. Assume that the prescribed monodromy representation (\ref{repre})  is irreducible.  Then the
generalized Riemann-Hilbert problem for $\mm$ has a 
solution if  one at least of the prescribed singularities is without roots.
\end{corollary}

\section{Reducible Solutions of the Generalized Riemann-Hilbert Problem} We  now consider  a non-generic datum $\mm_s$  for which the
generalized Riemann-Hilbert problem has a  solution.
One may  ask whether it is possible to realize this datum by a {\it reducible}  system
(\ref{sinv}) of differential equations of the form (\ref{red}). The
following statement  is a generalization of the main
result of \cite{Ma1} and  answers  the question.

\begin{theorem}    \label{mainred}
Let  $\mm_s$ be a non-generic datum. Then, under the conditions of 
Theorem 2, the reduced datum $\mm$ can be realized by a reducible system of the form (\ref{red}).

\end{theorem}
\begin{proof} We  proceed  as for the proof of Theorem 2, and keep the same notation as before. Consider  a vector bundle $(F^{\Lambda}, \nabla^{\Lambda})$ in $\ee$, and assume that
the pair $(F^{\Lambda}, \nabla^{\Lambda})$ is stable and holomorphically trivial, hence of degree zero.

\medskip
\noindent
The idea in the first step of the proof is the following. Starting with  $(F^{\Lambda}, \nabla^{\Lambda})$ we construct a pair
$(F^{ \tilde{\Lambda} }, \nabla^{ \tilde{\Lambda} })$ in $\ee$ which has a subbundle
$\tilde{F}_1$ stabilized by the connection $\nabla^{ \tilde{\Lambda} }$. The construction must be carried out in such a way that the pairs $( \tilde{F}_{1}, \nabla^{ \tilde{\Lambda} }|_{ \tilde{F}_{1} } )$ and $( F^{ \tilde{\Lambda} } / \tilde{F}_{1} , \nabla^{ \tilde{\Lambda} }_{q} )$,  where $\nabla^{ \tilde{\Lambda} }_{q}$ is the connection  induced  on the quotient bundle $F^{ \tilde{\Lambda} } / \tilde{F}_{1}$, are stable, that they have degree zero and that the difference between any two entries of the matrix $\tilde{ \Lambda }_1$ is greater than $(R-2+n)p$. As before, $R$ denotes the sum of all Poincar\'e ranks. For the construction of $\nabla^{ \tilde{\Lambda} }_{q}$, see the beginning of section 3.

\medskip
\noindent
The construction can be achieved as follows. Since the reduced datum $\mm_s$ is non-generic,  the set $\cal{F}$ of
proper subbundles of $F^{ \Lambda }$ that are stabilized by the connection $\nabla^{ \Lambda}$
is non-empty. Consider a bundle $F_1$ in $\cal{F}$ of maximal rank with the property that
$$ \mathrm{deg}(F_1) = \max_{ F' \in \cal{F} } \mathrm{deg}(F'). $$
It follows from the stability of $(F^{\Lambda}, \nabla^{\Lambda})$  that $\mathrm{deg}(F_1) < 0$.
Consider any proper filtration $F_{t_1} \subset F_{1} \subset F_{t_2} \subset F^{ \Lambda }$,
where $F_{t_1}$ and $F_{t_2}$ belong to $\cal{F}$ (by a {\it proper} filtration we mean that all inclusions are strict). In the following,
such a filtration will be called a {\it stabilized filtration}, and it satisfies the following inequalities:
$$ \mathrm{deg}( F_{t_1} ) \leq \mathrm{deg}(F_1) \ \ , \ \ 
\mathrm{deg}( F_{t_2} ) \leq \mathrm{deg}(F_1) - 1 \ \ , \ \ \mathrm{deg}( F^{ \Lambda } ) = 0. $$
Let $\lambda_{i}^{j}$ be an  entry of a block $\Lambda_{i}^{l}$ of $\Lambda_i$ 
corresponding to the superblock $\tilde{H}_{l}$ without roots of $\tilde{H}$ (see (\ref{adroots})). As in the
proof of Theorem 2, we have for any $N \in \mathbb{Z}$ 
$$ N( \lambda_{i}^{j} + \rho_{i}^{j} ) = \tilde{ \lambda }_{i}^{j} + \rho_{i}^{j} + \alpha_{i}^{j} $$
where $\tilde{ \lambda }_{i}^{j} \in \mathbb{Z}$ and $ 0 \leq \mathrm{Re}\  \alpha_{i}^{j} < 1$.
If the block $\Lambda_{i}^{l}$ corresponds to a factor  of dimension
$m_{il}$ of a superblock $\tilde{H}_{l}$
with roots, we replace this block by the block
${\Lambda'}_{i}^{l} = N \Lambda_{i}^{l} + s I_{ m_{il} }$ where
$$ s = \left[ \frac{ \sum \mathrm{Re} (N-1) \rho_{i}^{j} }{ m_{il} } \right] $$
and where the sum is taken over all eigenvalues of the superblock $\tilde{J}_{j}$ of $(\ref{adroots})$. To get
a bundle of degree  zero with  large enough differences between any two exponents, we
shall modify the integers $\tilde{ \lambda }_{i}^{m}$ as follows. We choose integers $k_{1}^{m}$,
 $1 \leq j \leq p$,  such that $k_{1}^{l_1} - k_{1}^{l_2} > (R-2+n)p$ \ if for $l_1 < l_2$ the entry $e_{l_1,l_2}$ of
$\tilde{J}$ is non-zero, and such that $k_{1}^{j} > 2np$\ for $1 \leq j \leq p$. We define
$$ \begin{array}{l}
             \displaystyle {\lambda'}_{1}^{1} = \tilde{ \lambda }_{1}^{1} + \sum_{i,j} \alpha_{i}^{j} +
r + k_{1}^{1} \\
             \displaystyle {\lambda'}_{1}^{j} = \tilde{ \lambda }_{1}^{j} + k_{1}^{j} \ \ , \ \
{\lambda'}_{1}^{p} = \tilde{ \lambda }_{1}^{p} - \sum_{j=1}^{p} k_{1}^{j} \\
             \displaystyle {\lambda'}_{i}^{j} = \tilde{ \lambda }_{i}^{j} \ \ \mbox{for all other $(i,j)$,}
                \end{array} $$
with
$$r = \sum m_{il} \left \{ \frac{ \sum \mathrm{Re} (N-1) \rho_{i}^{j} }{ m_{il} }  \right\},$$
where
the sum is taken over all superblocks $\tilde{J}_{l}$ with roots, for $2 \leq i \leq n$, and
$\{ \  \}$ stands for the fractional part. The matrices ${\Lambda'}_{i}$
with entries ${\lambda'}_{i}^{j}$ are clearly admissible,
with the additional property that the difference between any two diagonal entries of ${\Lambda'}_{1}$ is greater than $(R-2+n)p$.

\medskip
\noindent
Consider the pair $(F^{\Lambda'}, \nabla^{\Lambda'})$ obtained for some $N>>1$. For any stabilized 
filtration ${F}_{t_1}' \subset {F}_{1} '\subset {F}_{t_2}'\subset F^{ \Lambda' }$, we get
$$ \mathrm{deg}( {F}_{t_1}' ) < \mathrm{deg}( {F}_{1}' ) = N\  \mathrm{deg}( F_1 ) + c_{n,p,k} < 0,$$
$$  \mathrm{deg}( {F}_{t_2}' ) \leq N\  \mathrm{deg}( F_1 ) - N
+ 2np + \sum_{j=1}^{p} k_{j}^{1} < \mathrm{deg}( {F}_{1}') $$
and
$$ \mathrm{deg}( F^{ \Lambda' } ) = 0, $$
where $c_{n,p,k}$ is a sum of terms
involving $\alpha_{i}^{j}$,
$m_{il} \left\{ \sum \mathrm{Re} (N-1) \rho_{i}^{j} / m_{il} \right\}$ and $k_{1}^{j}$.
To construct the pair $(F^{ \tilde{\Lambda} }, \nabla^{ \tilde{\Lambda} } )$, we shall modify
the integers
${\lambda'}_{i}^{j}$ as follows
$$  \begin{array}{l}
          \tilde{ \lambda }_{1}^{1} = {\lambda'}_{1}^{1} + N\  \mathrm{deg}(F_1) + c_{n,p,k} \\
          \tilde{ \lambda }_{1}^{p} = {\lambda'}_{1}^{p} - N\  \mathrm{deg}(F_1) - c_{n,p,k} \\
          \tilde{ \lambda }_{i}^{j} = {\lambda'}_{i}^{j} \ \  \mbox{for the other $(i,j)$.}
                 \end{array}  $$
Again,  the matrices $\tilde{ \Lambda }_{i}$
with entries $\tilde{ \lambda }_{i}^{j}$ are admissible
with the additional property that the difference between any two diagonal entries of $\tilde{ \Lambda }_{1}$ is greater than $(R-2+n)p$. The pair $(F^{ \tilde{\Lambda} }, \nabla^{ \tilde{\Lambda} })$
moreover satisfies the following property. For any stabilized filtration
$\tilde{F}_{t_1} \subset \tilde{F}_{1} \subset \tilde{F}_{t_2} \subset F^{ \tilde{ \Lambda} }$ we have
$$ \mathrm{deg}( \tilde{F}_{t_1} ) < 0 \ \ , \ \ \mathrm{deg}( \tilde{F}_{1} ) = 0 \ \ , \ \
\mathrm{deg}( \tilde{F}_{t_2} ) < 0 \ \ , \ \ \mathrm{deg}( F^{ \tilde{ \Lambda} } ) = 0 . $$
This says that the pairs $( \tilde{F}_{1}, \nabla^{ \tilde{\Lambda} }|_{ \tilde{F}_{1} } )$ and 
$( F^{ \tilde{\Lambda} } / \tilde{F}_{1} , \nabla^{ \tilde{\Lambda} }_{q} )$ are stable and of degree zero. 

\medskip
\noindent In the second part of the proof, consider a meromorphic trivialization
of the bundle
$(F^{ \tilde{\Lambda} }, \nabla^{ \tilde{\Lambda} } )$ holomorphic outside of the point $a_{1}=0$, that induces a meromorphic trivialization of $\tilde{F}_{1}$.
The corresponding global system $(2)$ constructed from this trivialization has
the prescribed generalized monodromy data, except at $a_{1}=0$, where the Poincar\'e rank may be greater than $r_{1}$. We can choose a formal fundamental matrix $\hat{Y}_{1}(z)$ of $(2)$ of the form
$$ \hat{Y}_{1}(z) = \left( \begin{array}{cc}
                                                 \hat{Y}_{1}^{1}(z) & \ast \\
                                                           0 &  \hat{Y}_{2}^{1}(z)
                                             \end{array} \right) = \hat{F}(z) z^{ \tilde{ \Lambda }_{1} } z^{ \tilde{J} } e^{Q(z) },\ \ \  \hat{F} = \left( \begin{smallmatrix} \hat{F}_{1} & \ast \\ 0 & \hat{F}_{2} \end{smallmatrix} \right),$$
where  $ \hat{F}$   is  formal  meromorphic  and where the matrix $\hat{Y}_{1}^{1}(z)$ is chosen to be the formal fundamental matrix of a
subsystem $\omega_{1}$ that represents the restriction of $\nabla^{ \tilde{\Lambda} }$ on
$\tilde{F}_{1}$, and $\hat{Y}_{1}^{2}(z)$  the formal fundamental matrix of a quotient-system $\omega_{2}$ that represents the connection $\nabla^{ \tilde{\Lambda} }_{q}$ on
$F^{ \tilde{\Lambda} } / \tilde{F}_{1}$.

\medskip
\noindent
The bundles $\tilde{F}_{1}$ and $F^{ \tilde{\Lambda} } / \tilde{F}_{1}$ are holomorphically equivalent to sums
of line bundles

\begin{equation} \label{ast}
\tilde{F}_{1} \cong {\cal{O}} (c_{1}) \oplus \cdots \oplus {\cal{O}} ( c_{p_1} ) \ \ , \ \ 
      F^{ \tilde{\Lambda} } / \tilde{F}_{1} \cong {\cal{O}}( c_{p_{1} + 1} ) \oplus \cdots \oplus {\cal{O}}( c_p ) 
\end{equation}

\noindent where $c_{1} \geq \ldots \geq c_{p_1}$ and $c_{p_{1} + 1} \geq \ldots \geq c_p$. 

\medskip
\noindent
Since they are of degree zero, the trace of the matrices $K_{1} = \mathrm{diag}( c_{1},\ldots,c_{p_1} )$ and
$K_{2} = \mathrm{diag}( c_{p_{1} + 1},\ldots,c_{p} )$ is equal to zero. By construction, the bundles
$( \tilde{F}_{1}, \nabla^{ \tilde{\Lambda} }|_{ \tilde{F}_{1} } )$ and 
$( F^{ \tilde{\Lambda} } / \tilde{F}_{1} , \nabla^{ \tilde{\Lambda} }_{q} )$ are stable  pairs. By Lemma~2, we get the following estimates

\begin{equation} \label{est}
| c_{j} - c_{k} | \leq (R + n - 2)p
\end{equation} 
for $1 \leq j,k \leq p$,
It follows from the holomorphic equivalence  $( \ref{ast})$ that there is a gauge transformation
$\hat{Y}_{1}^{b}(z) = \Gamma^{b}(z)  \hat{Y}_{1}(z)$,  holomorphic outside of zero, of the form
$$ \Gamma^{b}(z) = \left( \begin{array}{cc}
                                                   \Gamma_{1}^{b}(z) & 0 \\
                                                                    0 & \Gamma_{2}^{b}(z)
                                               \end{array} \right), $$
that transforms the system $(2)$ into a system with the following formal fundamental solution at zero $$ \hat{Y}_{1}^{b}(z) = \left( \begin{array}{cc}
                                                 z^{-K_{1}} \hat{F}_{1}^{0}(z) & \ast \\
                                                           0 & z^{-K_{2}} \hat{F}_{2}^{0}(z) 
                                             \end{array} \right) z^{ \tilde{ \Lambda }_{1} } z^{ \tilde{J} } e^{Q(z) }, $$ 
where $\hat{F}_{j}^{0}$, $j=1,2$,  is  formal-holomorphically invertible. From Lemma 4  we get a gauge
transformation $ \hat{Y}_{1}^{k}(z) = \Gamma(z) \hat{Y}_{1}^{b}(z) $,  holomorphic outside of zero,
that transforms the latter system into a system with the following formal fundamental
solution at zero,
$$ \hat{Y}_{1}^{k}(z) =  \left( \begin{array}{cc}
                                                        \hat{H}_{1}(z) & \hat{H}_{3}(z) \\
                                                                       0 & \hat{H}_{2}(z)
                                                     \end{array} \right) z^{ \tilde{K} + \tilde{ \Lambda }_{1} } z^{ \tilde{J} } e^{Q(z) }, $$
where $\tilde{K}$ is obtained after a suitable permutation of the diagonal elements of
$\mathrm{diag}( K_{1} , K_{2} )$, and where $\hat{H}_{1}$, $\hat{H}_{2}$ are formal-holomorphically invertible, and $\hat{H}_{3}$ is formal meromorphic. Moreover, the estimates on the $c_{j}$ and
$\lambda_{1}^{j}$ imply that $\tilde{K} + \tilde{ \Lambda }_{1}$ is admissible.

\medskip
\noindent
We now need the following lemma, which is given in its analytic version in \cite{Bol2}.
\begin{lemma} Consider a formal meromorphic matrix $F(z)=\bigl({}^{  F_{1}(z)}_{    F_{2}(z)}\bigr)$
where $F_{2}( z )$ is  formal-holomorphically invertible. There exists a
meromorphic matrix $\Gamma(z)$ at 0, which is  holomorphically invertible outside of zero, such that
$$ \Gamma( z ) F( z ) = \left( \begin{array}{c}
                                    \tilde{F}_{1}( z ) \\
                                    F_{2}( z )
                              \end{array} \right) ,$$
where $\tilde{F}_{1}$ is formal holomorphic.
\end{lemma}
Using this lemma, we gauge-transform the system into a system $dy = \tilde{ \omega }y$ that has $\hat{Y}_{1}^{k}(z)$ as  formal fundamental solution, and 
where $\hat{H}_{3}$ is  formal holomorphic.

\medskip
\noindent
Formula  $(\ref{calc})$ shows that the latter system has Poincar\'e rank $r_1$ at zero, and the required Poincar\'e ranks at
all other points. Moreover,  the form of the fundamental matrix $\hat{Y}_{1}^{k}(z)$ tells us that this system
is \emph{reducible}, which means that the coefficient matrix $\tilde{ \omega }$ is  upper block-triangular,
$$ \tilde{ \omega } = \left( \begin{array}{cc}
                                                \tilde{ \omega }_{1} & \ast \\
                                                                0 & \tilde{ \omega }_{2}
                                                 \end{array} \right). $$

\end{proof}
As an application of the preceding results, we will consider the problem of reducibility for a
special type of systems $(2)$ which we will call \emph{formally fuchsian}.
\begin{definition} A differential system $(2)$ is called \emph{formally fuchsian} on
$\mathbb{P}^{1}( \mathbb{C} )$ if its formal
fundamental solution  $(\ref{Ffsol})$ at each singular point $a_{i}$, $1 \leq i \leq n$, is
$$ \hat{Y}_{i}(z) = \hat{F}_{0}(z) (z-a_{i})^{ \Lambda_{i} } (z-a_{i})^{ \tilde{J} } e^{ Q(z) } $$
where $\Lambda_{i}$ is admissible and $\hat{F}_{0}$ is formal-holomorphically invertible.
\end{definition}
The following statement is a generalization of the main results of \cite{Ma1},   \cite{Ma2}.
\begin{proposition}
Consider a differential system (\ref{sinv}) which is  formally fuchsian on
$\mathbb{P}^{1}( \mathbb{C} )$. Assume  that the generalized monodromy data of (\ref{sinv}) define a non-generic datum $\mm_s$
and that one at least of the singularities is without roots.
Then  the reduced datum $\mm$ 
can be realized by a reducible system of the form (\ref{red}).


\end{proposition}

\begin{proof} Consider the holomorphic trivial bundle $( F^{\Lambda}, \nabla^{ \Lambda } )$
constructed from the given system $(2)$. By construction, the degree $\mathrm{deg}( F^{ \Lambda } )$
of the bundle $F^{ \Lambda }$ is equal to zero. Moreover, for each subbundle $F'$ of
$F^{ \Lambda }$ that is stabilized by the connection $\nabla^{ \Lambda }$, the inequality
$\mathrm{deg}( F' ) \leq 0$ holds (see the beginning of section 3). If the pair $( F^{ \Lambda}, \nabla^{ \Lambda } )$ is stable,
then the result follows from Theorem 3. If the pair is unstable, there exists a proper subbundle $F_{1}$  stabilized by the
connection $\nabla$ and such that $\mathrm{deg}( F_1 ) = 0$. In section 3 we have seen that a holomorphic bundle on
$\pp$ is trivial if and only if it is a semi-stable bundle
of degree zero. From the holomorphic triviality of $F^{ \Lambda }$ we deduce that $F^{ \Lambda}$ is
semi-stable, which implies that $F_{1}$ and
$F^{\Lambda} / F_{1}$ are  semi-stable too,  of degree  zero. Therefore, $F_{1}$ and
$F^{\Lambda} / F_{1}$ are trivial bundles. To construct the reducible system we apply the
second part of the proof of Theorem 3 to the bundle
$( F^{\Lambda}, \nabla^{ \Lambda } )$. Indeed,  both bundles $F_{1}$ and $F^{\Lambda} / F_{1}$ have a trivial 
splitting type,  $c_{1} = \ldots = c_{p_1} = 0$ and $c_{p_{1} + 1} = \ldots = c_{p} = 0$,  hence the estimates (\ref{est}) hold.
\end{proof}

\section{The generalized Riemann-Hilbert problem in dimension two and three}

In this section, we solve the generalized Riemann-Hilbert problem in dimension two and three. We keep notation from previous sections. As usual, a formal solution (\ref{fsol}) is said to be convergent (divergent otherwise) if the asymptotic factor $\hat{F}$ in  (\ref{fsol}) is convergent in a neighbourhood of the singularity. We prove the following result.

\begin{theorem} \label{dt} In dimension two or three,  consider a reduced  datum $\mm_s$. The generalized Riemann-Hilbert problem  for $\mm$ has a  solution if we 
assume
that one at least of the prescribed singularities is without roots, and that the formal fundamental solution of one at least of the local systems $(\ref{si})$  is divergent.

\end{theorem}

\medskip
\noindent
{\bf Remark} If all the local systems $(\ref{si})$,  $1 \leq i \leq n$, have a convergent fundamental solution, then the generalized Riemann-Hilbert  problem  reduces to the classical Riemann-Hilbert problem. Indeed,  a finite number of gauge transformations of the form  $y={\mathrm e}^{q(z-a_i)}u$
(modulo a M\"obius transform if $a_i=\infty$) where $q(t)\in \frac{1}{t}\cc[\frac{1}{t}]$, will reduce the datum $\mm$ to
a datum of fuchsian singularities only. This is due to the fact that there is no Stokes phenomenon at the irregular singularities in this case, hence the exponential part $Q$  is a scalar matrix  at each $a_i$.
The classical problem always has a solution in dimension two (\cf \cite{AB}). In dimension three, a complete classification of the counterexamples for the classical problem was given in \cite{AB} and \cite{Glad}. Thus, the GRH-problem always has a solution in dimension two and is completely elucidated in dimension three, if we except the case where all data are those of irregular singularities with roots.

\medskip
\noindent
In dimension two the result is a generalization of Theorem
1 of \cite{JLP1}. In dimension three, it is a generalization of the main result of \cite{Ba2}.

\subsection{Proof of Theorem \ref{dt} in dimension two}

Choose a set
$\Lambda = ( \Lambda_{1}, \ldots, \Lambda_{n} )$ of admissible matrices and consider the
extension $(F^{\Lambda}, \nabla^{\Lambda})$ obtained {\it via}  the construction explained earlier . There are three parts in the proof.

\medskip \noindent
We first assume that the bundle $F^{\Lambda}$ has no proper subbundle  stabilized
by the connection $\nabla^{\Lambda}$. Thus, the pair $(F^{\Lambda}, \nabla^{\Lambda})$ is stable and by Theorem 2 the GRH-problem for $\mm$ has a solution. 

\medskip
\noindent
Now  assume that the bundle
$F^{\Lambda}$ is the direct sum of two proper subbundles $F_{1}$ and $F_{2}$, that are stabilized by $\nabla^{\Lambda}$. Starting with the pairs $( F_{k}, \nabla^{\Lambda}|_{F_k} )$, $k=1,2$, one can easily construct two global differential systems
$$ dy = \omega_{1}y \ \ , \ \ dy = \omega_{2}y $$
on $\pp$ with the generalized monodromy data of the bundles
$( F_{1}, \nabla^{\Lambda}|_{F_1} )$ and $( F_{2}, \nabla^{\Lambda}|_{F_2} )$ respectively. It is
easy then to see that the differential system
$$ dy = \left( \begin{array}{cc}
                                 \omega_{1} & 0 \\
                                            0 & \omega_{2}
                         \end{array} \right) y$$
has the prescribed generalized monodromy data $\mm$.

\medskip
\noindent
In the last part of the proof we assume that the bundle $F^{\Lambda}$ has a unique proper subbundle $F_{1}$ that
is stabilized by  $\nabla^{\Lambda}$.
As in  $(\ref{redomega})$, $(\ref{redfsol})$, starting with the local differential system $dy = \omega^{\Lambda_{i_1}}y$, we  construct
a local differential system $dy = \omega_{i_1}'y$ with a formal fundamental solution of the form
$$ \hat{Y}_{i_1}(z) = \left( \begin{array}{cc}
                                                  \hat{Y}_{i_1}^{1} & \ast \\
                                                           0 &   \hat{Y}_{i_1}^{2}
                                               \end{array} \right) =
\hat{F}_{0}(z) ( z-a_{i_1} )^{ \Lambda_{i_1} } ( z-a_{i_1} )^{J'} e^{Q'(z)}, $$
where $ \hat{Y}_{i_1}^{1}$ is a formal fundamental solution of a local system defining the restriction of
the connection $\nabla^{\Lambda}$ to the subbundle $F_{1}$. {\it Via} a basis change
$\hat{Y}_{i_1}S$,  $S \in \mathrm{GL}(2, \mathbb{C})$, and a suitable permutation of the
diagonal elements of $\Lambda_{i_1}$ (preserving admissibility) we may assume that $J'$ has
two specified forms which we detail below.

\medskip
\noindent
{\it Case 1.}  The matrix $J'$ has the form

$$ J' = \left( \begin{array}{cc}
                            \rho_{i_1}^{1} & 1 \\
                                 0 & \rho_{i_1}^{2}
                       \end{array} \right) . $$
In this case, the matrix $Q'(z)$ is a scalar matrix of the form $Q'(z) = q'(z)I_{2}$. From the classical theory
(\cf \cite{BJL2},  p.262), we know that the gauge transformation $u = \exp(-q'(z))I_{2}y$ changes
the system $dy = \omega_{i_1}'y$ into a local system with a regular singularity at $a_{i_1}$. This contradicts our assumptions, hence case 1 does not occur.

\medskip
\noindent 
{\it Case 2.}  The matrix $J'$ is of the form

$$ J' = \left( \begin{array}{cc}
                            \rho_{i_1}^{1} & 0  \\
                                 0 & \rho_{i_1}^{2}
                       \end{array} \right) . $$
In this case, it is possible to construct a
stable pair $(F^{ \tilde{\Lambda} }, \nabla^{ \tilde{\Lambda}} )$ with the prescribed generalized
monodromy data. Thus, by Theorem 2, the GRH-problem is solved. Indeed, starting with the set of admissible matrices $\Lambda$,
we construct a new set of admissible matrices
$\tilde{\Lambda} = ( \tilde{\Lambda}_{1}, \ldots, \tilde{\Lambda}_{n} )$ in the following way. Let $b$ be a
positive integer. Let $\tilde{ \lambda }_{i_1}^{1} = \lambda_{i_1}^{1} - b$,
$\tilde{ \lambda }_{i_1}^{2} = \lambda_{i_1}^{2} + b$, and $\tilde{ \lambda }_{i}^{j} = \lambda_{i}^{j}$ for all
other $i,j$. For  sufficiently large $b$ we get the
inequality
$$ \sum_{i=1}^{n}  \tilde{\lambda}_{i}^{1} + \rho_{i}^{1} < \frac{1}{2} 
     \sum_{i=1}^{n} \sum_{j=1}^{2} \tilde{\lambda}_{i}^{j} + \rho_{i}^{j} =  \frac{1}{2} 
     \sum_{i=1}^{n} \sum_{j=1}^{2} {\lambda}_{i}^{j} + \rho_{i}^{j}. $$
Let $\tilde{\Lambda}_{i} = \mathrm{diag}( \tilde{\lambda}_{i}^{1}, \tilde{\lambda}_{i}^{2} )$, for $1 \leq i \leq n$, and consider the extension $( F^{\tilde{\Lambda}}, \nabla^{\tilde{\Lambda}} )$. The above inequality implies that the latter pair is stable.

\subsection{Proof of Theorem \ref{dt} in dimension three}

The proof in dimension three follows the same lines as in dimension two. We will only give the last part of the proof. 

\medskip
\noindent
We first assume that the bundle $F^{\Lambda}$ has a unique proper subbundle $F_{1}$ that
is stabilized by  $\nabla^{\Lambda}$. As in dimension two, and with the same notation, we shall construct a stable pair $(F^{ \tilde{\Lambda} }, \nabla^{ \tilde{\Lambda}} )$ with the given data, hence solve the GRH-problem in each case. To construct this pair we consider the following cases.

\medskip
\noindent {\it Case 1. } The matrix $J'$ has the form

$$ J' = \left( \begin{array}{ccc}
                            \rho_{i_1}^{1} & 0 & 0 \\
                                 0 & \rho_{i_1}^{2} & 1 \\
                                 0 & 0 & \rho_{i_1}^{3}
                       \end{array} \right)  \ \ \mbox{or} \ \  J' = \left( \begin{array}{ccc}
                            \rho_{i_1}^{1} & 0 & 0 \\
                                 0 & \rho_{i_1}^{2} & 0 \\
                                 0 & 0 & \rho_{i_1}^{3}
                       \end{array} \right).$$
Starting with  the set of admissible matrices $\Lambda$ we construct a set of admissible matrices
$\tilde{\Lambda} = ( \tilde{\Lambda}_{1}, \ldots, \tilde{\Lambda}_{n} )$ in the
following way. Let $b$ be a positive integer. Let $\tilde{\lambda}_{i_1}^{1} = \lambda_{i_1}^{1} - 2b$,
$\tilde{\lambda}_{i_1}^{2} = \lambda_{i_1}^{2} + b$, $\tilde{\lambda}_{i_1}^{3} = \lambda_{i_1}^{3} + b$ and
$\tilde{\lambda}_{i}^{j} = \lambda_{i}^{j}$ for all other $i,j$. For sufficiently large $b$ we get the following inequality
$$ \sum_{i=1}^{n} \sum_{j=1}^{2} \tilde{\lambda}_{i}^{j} + \rho_{i}^{j} < \frac{2}{3} 
     \sum_{i=1}^{n} \sum_{j=1}^{3} \tilde{\lambda}_{i}^{j} + \rho_{i}^{j} =  \frac{2}{3} 
     \sum_{i=1}^{n} \sum_{j=1}^{3} {\lambda}_{i}^{j} + \rho_{i}^{j}  .$$
Let $\tilde{\Lambda}_{i} = \mathrm{diag}( \tilde{\lambda}_{i}^{1}, \tilde{\lambda}_{i}^{2},
\tilde{\lambda}_{i}^{3} )$ for  $1 \leq i \leq n$, and
consider the extension $( F^{\tilde{\Lambda}}, \nabla^{\tilde{\Lambda}} )$. It follows from the above inequality  that this pair is stable, hence, by Theorem 2, the GRH-problem is solved in this case.

\medskip
\noindent {\it Case 2.}   The matrix $J'$ has the form

$$ J' = U^{-1} \left( \begin{array}{ccc}
                            \rho_{i_1}^{1} & 1 & 0 \\
                                     0 & \rho_{i_1}^{2} & 0 \\
                                      0 & 0 & \rho_{i_1}^{3}
                       \end{array} \right) U $$
where $U$ is the matrix described in section 1.2. From the set of admissible matrices $\Lambda$
we again construct a set of admissible matrices
$\tilde{\Lambda} = ( \tilde{\Lambda}_{1}, \ldots, \tilde{\Lambda}_{n} )$ in the
following way. Let $b$ be a positive integer, and let $\tilde{\lambda}_{i_1}^{1} = \lambda_{i_1}^{1} - b$,
$\tilde{\lambda}_{i_1}^{2} = \lambda_{i_1}^{2} - b$, $\tilde{\lambda}_{i_1}^{3} = \lambda_{i_1}^{3} + 2b$ and
$\tilde{\lambda}_{i}^{j} = \lambda_{i}^{j}$ for all other $i,j$. For sufficiently large $b$ we get the inequality
$$ \sum_{i=1}^{n} \sum_{j=1}^{2} \tilde{\lambda}_{i}^{j} + \rho_{i}^{j} < \frac{2}{3} 
     \sum_{i=1}^{n} \sum_{j=1}^{3} \tilde{\lambda}_{i}^{j} + \rho_{i}^{j} =  \frac{2}{3} 
     \sum_{i=1}^{n} \sum_{j=1}^{3} {\lambda}_{i}^{j} + \rho_{i}^{j} $$
Consider the extension $( F^{\tilde{\Lambda}}, \nabla^{\tilde{\Lambda}} )$ obtained with $\tilde{\Lambda}_{i} = \mathrm{diag}( \tilde{\lambda}_{i}^{1}, \tilde{\lambda}_{i}^{2},
\tilde{\lambda}_{i}^{3} )$ for $1 \leq i \leq n$. It is a stable pair in view of the above inequality. Thus, by Theorem 2,  the GRH-problem has a solution in this case.

\medskip
\noindent 
{\it Case 3. } The matrix $J'$ is a Jordan block

$$ J' = \left( \begin{array}{ccc}
                            \rho_{i_1}^{1} & 1 & 0 \\
                                     0 & \rho_{i_1}^{2} & 1 \\
                                      0 & 0 & \rho_{i_1}^{3}
                       \end{array} \right) . $$
In this case, the matrix $Q'(z)$ is a scalar matrix of the form $Q'(z) = q'(z)I_{3}$s and the gauge transformation $u = \exp(-q'(z))I_{3}y$ changes the system
$dy = \omega_{i_1}'y$ into a system with a regular singularity at $a_{i_1}$. This
contradicts our assumptions,  hence this case will not occur. 

\medskip
\noindent
Now ssume that the bundle $F^{\Lambda}$ has a stabilized filtration $F_{1} \subset F_{2} \subset F^{\Lambda}$ where $F_{1}$ has dimension 1 and $F_{2}$  dimension 2. As in $(\ref{redomega})$, $(\ref{redfsol})$, starting with the local differential system $dy = \omega^{ \Lambda_{i} }y$ we construct a local differential system $dy = \omega_{i}'y$ which has a formal fundamental solution of the form
$$ \hat{Y}_{i}(z) = \left( \begin{array}{ccc}
                                                  \hat{Y}_{i}^{1} & \ast & \ast \\
                                                           0 &   \hat{Y}_{i}^{2} & \ast \\
                                                            0 & 0 & \hat{Y}_{i}^{3}
                                               \end{array} \right) =
\hat{F}_{0}(z) ( z-a_{i} )^{ \Lambda_{i} } ( z-a_{i} )^{ J'_{i} } e^{Q'(z)} $$
where $ \hat{Y}_{i}^{1}$ is a formal fundamental solution of a local system which defines the restriction of
the connection $\nabla^{\Lambda}$ to the subbundle $F_{1}$, and $\hat{Y}_{i}^{2}$ is a formal
solution of a local system which defines the connection constructed from $\nabla^{ \Lambda }$
on the quotient bundle $F_{2} / F_{1}$. We notice that the matrices $J'_{i}$, $1 \leq i \leq n$ are
simultaneously upper-triangular. 

\medskip
\noindent
The following lemma can be proved in the same way as Proposition 3.1.3 of \cite{Van}.
\begin{lemma}
There exist invertible upper triangular matrices $S_{i}$, $1 \leq i \leq n$, and integers
$\varphi_{i}^{j}$, $1 \leq i \leq n$, $j=1,2,3$,  such that  \\
a) if, for $l_1 < l_2$, the entry $e_{l_1,l_2}$ of $S_{i}^{-1} J'_{i} S_{i}$  is non-zero,
then $\varphi_{i}^{l_1} \geq \varphi_{i}^{l_2}$,\\
b) for $t=1,2,3$, one has
$$ \sum_{i=1}^{n} \varphi_{i}^{t} + \rho_{i}^{t} = 0. $$
\end{lemma}
In the following, let $\tilde{ \Lambda }_{i} = \mathrm{diag}( \varphi_{i}^{1}, \varphi_{i}^{2},
\varphi_{i}^{3} )$ and  consider the pair $(F^{ \tilde{ \Lambda } }, \nabla^{ \tilde{ \Lambda } } )$.
There is a stabilized filtration
$\tilde{F}_{1} \subset \tilde{F}_{2} \subset F^{ \tilde{ \Lambda } }$ such that
$\mathrm{rank}( \tilde{F}_{1} ) = 1$ and  $\mathrm{rank}( \tilde{F}_{2} ) = 2$. The bundles
$\tilde{F}_{1}$, $\tilde{F}_{2} / \tilde{F}_{1}$, $F^{ \tilde{ \Lambda } } / \tilde{F}_{2}$ are equivalent
to line bundles,
$$ \tilde{F}_{1} \cong {\cal{O}} (c_{1}) \ \ , \ \  \tilde{F}_{2} / \tilde{F}_{1} \cong {\cal{O}} (c_{2}) \ \ , \ \
F^{ \tilde{ \Lambda } } / \tilde{F}_{2} \cong {\cal{O}} (c_{3}). $$
By the above lemma, these bundles are of degree zero, that is, 
$c_{1} = c_{2} = c_{3} = 0$. To construct a global system on  $\pp$
with the given  data, we  use the second part of the proof of
Theorem 3 applied to the pair $(F^{ \tilde{ \Lambda } }, \nabla^{ \tilde{ \Lambda } } )$. Indeed, each of the bundles
$\tilde{F}_{1}$, $\tilde{F}_{2} / \tilde{F}_{1}$, $F^{ \tilde{ \Lambda } } / \tilde{F}_{2}$ has a trivial
splitting type, which implies that all the estimates (\ref{est})  hold. This ends the proof of Theorem \ref{dt} in dimension three.

\section{The global inverse problem in differential Galois theory}

In this section we  show how the global inverse problem in differential Galois theory is related to the generalized Riemann-Hilbert problem. We first recall some  results of differential Galois theory, and we refer to the book of M. Singer and M. van der Put \cite{PS} for an extensive exposition of the theory.

\subsection{Differential Galois groups}
Consider a linear differential system 
\begin{equation}   \label{gs}
y'=By 
\end{equation}
where $B$ is a $p\times p$ matrix with entries in a differential field $K$, whose subfield $C$ of constants is algebraically closed and of characteristic zero. A {\it Picard-Vessiot extension} of $K$ with respect to (\ref{gs}) is a differential extension of $K$ with no new constants, containing the entries of a fundamental solution of (\ref{gs}) and generated by these entries over $K$. Such extensions always exist, and they are isomorphic. The {\it Galois group} of (\ref{gs})  over $K$ is the group $G$ of all differential $K$-automorphisms of a Picard-Vessiot extension of $K$ with respect to (\ref{gs}). A representation of this group in $\gl(p,C)$ is given by any fundamental solution of (\ref{gs})  generating a Picard-Vessiot extension of $K$ for (\ref{gs}), and $G$ is then an algebraic subgroup of $\gl(p,C)$. Systems which are $K$-equivalent have isomorphic Galois groups over $K$.

\medskip
\noindent
The  inverse problem in differential Galois theory is the following: {\it Given a differential field $K$ as before, and a linear algebraic group $G$ defined over the field $C$ of  constants of $K$, is it possible to realize $G$ as a differential Galois group over $K$? }

\medskip
\noindent
Over the field $C(z)$ of rational functions  or the field $\cc(\{z\})$ of convergent Laurent series, the problem is completely  solved  (\cf \cite{TT}, \cite{R4}, \cite{R5}, \cite{R6}, \cite{JH}, see also \cite{Kov1}, \cite{Kov2},  \cite{Sin}, \cite{MS1}, \cite{MS2}, \cite{MS3}, \cite{Be}, \cite{CMS}, and \cite{PS} chap. 11). Any group can be realized as a differential Galois group over $C(z)$, but not necessarily as a Galois group over  $\cc(\{z\})$ as we will see now. 

\subsection{The local inverse problem}
In this section, the differential base field is $\cc(\{z\})$.
\subsubsection{Ramis's solution}
The  inverse problem over $\cc(\{z\})$, also called  the {\it local Galois inverse problem}, was solved by J.-P. Ramis (\cite{R3},  \cite{R4}, \cite{R5}, \cite{R6}, \cf  \cite{MS2},  see also \cite{PS},  chap.11), who has proved that a linear algebraic group $G$ is a Galois group over $\cc(\{z\})$  if and only if it has a {\it local Galois structure}, which he defined as  follows. It is a triple ${\mathcal L}=(T,a,{\mathcal N})$, where
 
\noindent
(i) \ $T$ is a torus of $G$ and $a\in G$  normalizes  $T$

\noindent
(ii) \ the image of $a$  generates the finite group $G/G^0$ (with the usual notation $G^0$ for the identity component of $G$)

\noindent
(iii) \ ${\mathcal N}$  is a  Lie subalgebra of dimension $\le 1$ of the Lie algebra $\mathcal G$ of $G$, which commutes with $a$ and with $T$

\noindent
(iv) \ $\mathcal G={\cal T}+{\cal N}+{\cal Q}(T)$, where $\cal T$ denotes the Lie algebra of $T$, and ${\cal Q}(T)$  the critical subalgebra  for $T$, defined as the Lie subalgebra of $\cal G$ generated  by the rootspaces of $\cal G$ under the adjoint action of $T$.
 
\medskip
\noindent
In \cite{MS2} it was shown that local Galois groups  also are  characterized by the condition (i)  $G/G^0$ is cyclic, (ii) ${\mathrm {dim}}(R_u/(R_u,G^0))\le 0$, and (iii) $G/G^0$ acts trivially on $R_u/(R_u,G^0)$, where $R_u$ is the unipotent radical of $G$.

\medskip
\noindent
A {\it reduced} local Galois structure on $G$ is a local Galois structure ${\mathcal L'}=(T',a',{\mathcal N}')$ as before, satisfying the additional conditions  that

\noindent
{\it (i)} $T'$ is a maximal torus

\noindent
{\it (ii)}  $a'$  is of finite order in $G$ (and  semisimple)

\smallskip
\noindent
{\it (iii)} ${\mathcal N}'$  is either $(0)$ or the Lie algebra of a subgroup isomorphic to $\cc$, and ${\mathcal N}'\cap({\cal T}'+{\cal Q}(T'))=\emptyset$.

\medskip
\noindent
The density theorem of Ramis (\cite{R1}, \cite{R2}, \cf \cite{LR}, see also \cite{PS} theorem 11.13) states that the differential Galois group $G$ of a linear differential system (\ref{gs}) over $\cc(\{z\})$ is topologically generated by the formal monodromy, the Stokes matrices (as defined in section 2) and the {\it exponential torus}, that is, the torus $T_e$ of $K$-differential automorphisms of the field $K(e^{q_1},\ldots,e^{q_p})$  where the $q_i$'s (see section 1) are the diagonal entries of the exponential part $Q$ in (\ref{fsol}). 

\medskip
\noindent
The system (\ref{gs}) {\it via} its local Stokes data described in section 2, gives rise to a local Galois structure on $G$. In short, $T$ is the sum of  the exponential torus $T_e$ and the monodromy torus $T_m$ (generated by the semisimple part of the formal monodromy)  whereas the infinitesimal Stokes matrices (inverse images by the exponential map of the Stokes matrices) can be developed in $\cal G$ as sums of rootspace elements under the action of the exponential torus to produce generators of the critical subalgebra ${\cal Q}(T_e)$. The Lie algebra ${\mathcal N}$ arises from the unipotent part of the formal monodromy and $a$ from its finite part (\cf \cite{R3}, \cite{PS}).

\medskip
\noindent
Given any local Galois structure ${\mathcal L}=(T,a,{\mathcal N})$ on $G$, in particular one induced  by a system (\ref{gs}), there is a natural way to associate to $\ll$ a reduced local Galois structure ${\mathcal L'}=(T',a',{\mathcal N}')$, where $T'$ contains $T$, and $a'$ equals $a$ modulo $G^0$.

\medskip
\noindent
\subsubsection{The Poincar\'e rank for a local Galois structure}
Our aim here is  to realize local Galois data with a minimal Poincar\'e rank.
\begin{definition}
A  {\rm local Galois datum} is a pair $(G,\ll)$ where $G$ is a linear algebraic group endowed with a reduced local Galois structure $\ll$. A differential system (\ref{gs}) over $\cc(\{z\})$ is said to {\rm realize} $(G,\ll)$ if $G$ is the differential Galois group of (\ref{gs}) over $\cc(\{z\})$ and if the system (\ref{gs}) induces the local Galois structure $\ll$ on $G$.

\noindent
The {\rm Poincar\'e rank} $r_{\ll}$ of a local Galois datum $(G,\ll)$ is the smallest possible Poincar\'e rank of a differential system realizing these data. For a given group $G$ with local Galois structures,  let  $r(G)$ denote the minimal Poincar\'e rank $r_\ll$, over all possible local data $(G,\ll)$.
\end{definition}

\medskip
\noindent
To determine the Poincar\'e rank $r_{\ll}$ of a  given local datum $(G,\ll)$ we will use the construction of Ramis in his proof of the local inverse problem (\cite{R4}  section 2.1, \cf \cite{PS} pp. 273-74, 279-82).

\medskip
\noindent
We will carry out this construction  in such a way that the Katz rank $\rho_{\ll}$ of  the system (\ref{gs}) realizing the data, that is, the fractional degree in $z$ of the exponential part $Q$ (of a formal fundamental solution of (\ref{gs}) of the form (\ref{fsol}), section 1.2) is minimal. We determine $\rho_{\ll}$, and hence $r_{\ll}$, explicitely in terms of $\ll$.

\medskip
\noindent
Let $\ll=(T,a,{\cal N})$ be the given reduced local structure on $G$, where $a\in G$ acts by conjugation on the maximal torus  $T$ as an automorphim of order $\nu\in{\bf N}^*$. Let $G$ be given with a faithful representation $G\subset \gl(n,\cc)$ such that $T$ is a diagonal subgroup,  and let $\chi_i$, $i=1,\ldots,s,$ denote the corresponding distinct diagonal weights of $T$, which generate the (abelian)  dual group $\tt$ of $T$ as a $\zz$-module.

\noindent
Consider the $\qq$-vector space $E=\tt\oplus_{\zz}\qq$, and  the $\qq$-automorphism $\delta$ of $E$ of order $\nu$ induced by the conjugation by $a$ on $T$. The decomposition $\delta^{\nu}-id=\prod_{\nu' | \nu}\Phi_{\nu'}(\delta)$ where $\Phi_{\nu'}$ denotes the $\nu'$-th cyclotomic polynomial, yields a decomposition $E={\bigoplus_{k=1}^m}E_k$ of $E$ into a direct sum of $\delta$-invariant $\qq$-subspaces $E_k$, each  of dimension $\nu_k$  for some divisor $\nu_k$ of $\nu$, and such that $\Phi_{\nu_k}$ is the minimal polynomial of $\delta$ on $E_k$. 

\smallskip
\noindent
Let ${\cal F}=\bigoplus_{\lambda \in \qq,\lambda < 0} {\mathbb C}[z^{\lambda}]$ denote the ring of polynomials  in (non-negative) fractional powers of $1/ z$. For each $k$ one can realize $E_k$ as the $\qq$-span of an isomorphic image of some lattice $\bigoplus_{i=1}^{\varphi(\nu_k)}\zz p_i$, $p_i\in{\cal F}$, in the following way. For any given arbitrary integer $\mu_k\ge 1$ prime to $\nu_k$  we can choose $p_1=z^{-{\mu_k} /{\nu_k}}$ and $p_j=m^{j-1}(p_1)$, $j=1,\ldots,\varphi(\nu_k)$, where $m$ denotes the monodromy operator on $\cal F$. This defines an  isomorphism which clearly commutes with $\delta$ and $m$. The family $\underline p$ of all such polynomials $p_i\in\cal F$ for all $k$, is $m$-invariant and $\zz$-independent, and the above isomorphisms glue together in a global isomorphism $\psi : \cal P\rightarrow E$   from the $\qq$-span  $\cal P$ of $\underline p$ to $E$. Let  $c_i \in \cal P$, for each $i=1,\ldots,s$, denote the inverse image of $\chi_i$.

\noindent
The roots of the adjoint action of $T$ on $\cal G$ are elements of $\check T$ and each non-zero root is actually of the form $\chi_i-\chi_j$ for some $i,j \in \{1,\ldots,s\}, i\ne j$ (\cite{PS}, p.280, proof of lemma 11.16). For each non-zero root $\alpha$, let ${\cal G}_{\alpha}$ denote the corresponding rootspace.

\noindent
Note that since $\mu_k$ and  $\nu_k$ are relatively prime for all $k$, the correspondence $c_i \leftrightarrow \chi_i$ does not {\it a priori} depend on a precise choice of the $\mu_k$. For each $k$ let now  $\mu_k\ge 1$  be the smallest integer prime to $\nu_k$, and such that $\mu_k\ge {\rm dim}({\cal G}_{\alpha})$ for all $\alpha=\chi_i -\chi_j$ such that  the corresponding polynomial $c_i-c_j$ is of degree ${\mu_k}/{\nu_k}$.

\noindent
Let $\rho_{\ll}$ denote the largest of the fractional degrees $\mu_k/\nu_k$, for all divisors $\nu_k$ of $\nu$ occurring in the decomposition  of $E$.

\smallskip
\noindent
We have obtained the following result.

\begin{proposition} \label{katz}
Any local datum $(G,\ll)$ can be realized by a system (\ref{gs})  whose Katz degree is equal to
$\rho_\ll$.
\end{proposition}

\begin{proof}\  Let $(q_1,\ldots,q_n)$ denote the family of polynomials $q_i\in\cal F$ corresponding to the complete family of diagonal weights of $T$, taking into account their multiplicity and ordering in the given representation. In the construction of Ramis, this family  produces the exponential part $Q={\mathrm{ diag}}(q_1,\ldots,q_n)$ of the desired system, of degree ${\mathrm max}_k\{{{\mu_k}/{\nu_k}}\}$. The above choice of the $\mu_k$ makes it possible, in this construction, to define sufficiently  many Stokes operators to generate the critical algebra $\cal Q(T)$, hence to solve the local inverse problem with a minimal Katz degree. (The degree of any polynomial $q_i-q_j$ corresponding to a given root $\alpha$ is large enough to define sufficiently many Stokes rays).
\end{proof}

\medskip
\noindent
Conversely, any system (\ref{gs}) inducing local Stokes data $(G,\ll)$ can easily be seen to have a Katz degree greater or equal to $\rho_\ll$. In view of the previous construction, this implies the following result.

\begin{corollary} \label{minrk} 
With notation as above, we have $r_\ll=-[-\rho_\ll]$.
\end{corollary}

\subsection{The global inverse problem}

In this section we apply our results on the generalized Riemann-Hilbert problem to solve the global differential Galois problem with a better control of the singularities.

\medskip
\noindent
Let $G$ be a given linear algebraic group over $\cc$. We know that $G$ is the differential Galois group of some linear differential system over $\cc(z)$. Moreover, results of J.-P. Ramis also tell us that $G$ can be realized as the Galois group of systems with a certain type and number of singularities. 

\medskip
\noindent
Our aim is to prove the existence, under certain conditions, of  a system (realizing $G$ as its Galois group) with a minimal number of singularities, and with the smallest possible Poincar\'e rank at these.


\subsubsection{Number of singularities}
We first recall the main results of J.-P. Ramis about the singularities of a system solution of the global inverse problem, for which we refer to (\cite{R3}, Theorems 7.1.4 and 7.2.11, Proposition 7.1.6) and (\cite{PS}, Theorem 11.21).

\medskip
\noindent
Let $G$ be a complex algebraic group. Consider the subgroup $L(G)$ of $G$  generated by all (maximal) tori of $G$ and let $V(G)$ denote the quotient $G/L(G)$. Let $s=s(G)$ (resp. $\overline{s}=\overline{s}(G)$) denote the least positive integer $\ge 2$ such that $G$ (resp. $V(G)$) can be topologically genetated by $s-1$ (resp. $\sb -1$)  elements. By a {\it family of generators} of $G$ we  mean a family of  elements of $G$ generating $G$ topologically, and by a {\it  minimal} family of generators, a family of $s(G)-1$ such elements. By {\it generating} we will always, if not specified otherwise, mean ``generating topologically".

\medskip
\noindent
Tretkoff's original result, based on the Riemann-Hilbert correspondence, states that
any  linear algebraic group $G$ over $\cc$ can be realized as the Galois group over $\cc(z)$ of a differential system with $s(G)$ possible singularities, all  regular singular, and  fuchsian but possibly one. This was generalized by Ramis who proved that  a given linear algebraic group $G$  topologically generated by closed subgroups $G_1,\ldots,G_{m-1}$, $m\ge 2$, each endowed with a local Galois structure, is the Galois group over $\cc(z)$ of a system (\ref{gs}) with no more than $m$ singularities; these belong to a subset  $\{a_1,\ldots,a_m\}$ of $\pp$ such that the local Galois group of  (\ref{gs}) at each $a_i$ is $G_i$, $i=1,\ldots,m-1$, and $a_m$, if singular, is fuchsian. Moreover there exists such  a system for which all but one of its regular singularities are fuchsian. A more precise result states that any group $G$ is the Galois group of a differential system with $\sb(G)$ possible singularities, all fuchsian but one, possibly irregular. But we know nothing {\it a priori} sabout the  Poincar\'e rank at the irregular singularity.

\medskip
\noindent

\subsubsection{ Poincar\'e ranks}
In this section we apply our results on the generalized Riemann-Hilbert problem to refine the  results of Ramis, taking into account the Poincar\'e rank at the singularities.  

\medskip
\noindent
Throughout this section, $G$ is given with a faithful representation $G\subset\gl(n,\cc)$. We will use the following notation. For any family  $\mb=(M_1,\ldots,M_r)$ of elements of $G$, let $\chi_{\mb}$ denote the representation
$$ \chi_{\mb}:\pi_1(\mathbb{P}^1(\mathbb{C})\setminus\mathcal{D},z_0)\longrightarrow {\mathrm {GL}}(n,\mathbb C) $$
of the fundamental group of  ${\mathbb{P}}^1({\mathbb{C}})$ punctured at a set $\mathcal{D}=\{a_1,\ldots,a_{r+1}\}$ of  arbitrarily chosen   points of $\pp$ not containing the base-point $z_0$, such that  $\chi_{\mb}(\gamma_i)=M_i$, for all $i=1,\ldots,r$, where $\gamma_i$  denotes the class of an elementary loop around $a_i$.

\medskip
\noindent
The following result immediately follows from the Bolibrukh-Kostov irreducibility condition on the monodromy for the classical Riemann-Hilbert problem.
\begin{proposition} \label{galoisfuchs}
If the representation $G\subset\gl(n,\cc)$ is irreducible, then $G$ is the Galois group of a differential system with no more than $s(G)$ singularities, all fuchsian.
\end{proposition}
We will now combine the above mentionned results of differential Galois theory with results of the previous sections to get more precise results on the singularities of a system realizing a given group $G$ as its Galois group.

\begin{theorem} \label{galoisrecol}  Let $G$ be a linear algebraic group over $\cc$,  topologically generated by closed subgroups  $G_1,\ldots,G_{m-1}$, $m\ge 2$, each endowed with a reduced local Galois structure $\ll_i$, and let $\cd=\{a_1,\ldots,a_m\}$  be  an arbitrary set of $m$  points of $\pp$.  Let $\mm_s$ be a reduced datum  on $\cd$ realized by $m$  local systems  (\ref{si}) with Galois group $G_i$ and Poincar\'e rank $r_i$ at $a_i$, $i=1,\ldots,m-1$, and fuchsian at $a_m$. If $\mm$ fulfills the conditions of Theorem  \ref{main}, in particular if $\mm_s$ is generic, then $G$ is the Galois group of a global system with Poincar\'e rank $r_i$ at each $a_i$, $i=1,\ldots,m-1$, and which is fuchsian at $a_m$. 
\end{theorem}

\begin{proof} The datum $\mm$ is fuchsian  at $a_m$, hence without roots. Apply Theorem \ref{main}. A system (\ref{gs}) with  generalized monodromy data $\mm$ has $G_i$ as its local differential Galois group  at  $a_i$, $i=1,\ldots,m-1$, and since the $G_i$ together generate $G$, the global Galois group of (\ref{gs})  is $G$. If $\mm_s$ is generic, apply Corollary \ref{generic}.
\end{proof}



\noindent From Theorem \ref{galoisrecol} and from the monodromy criterion  of Corollary \ref{irr} we deduce the following results.
\begin{corollary} \label{galoisirr}  Let $G$ be a linear algebraic group over $\cc$,  topologically generated by closed subgroups  $G_1,\ldots,G_{m-1}$, $m\ge 2$, each endowed with a reduced local Galois structure $\ll_i$, and let $\cd=\{a_1,\ldots,a_m\}$  be  an arbitrary set of $m$  points of $\pp$. If there exists a family $\mb=(M_1,\ldots,M_{m-1})$ of elements of $G$ such that

\noindent {\it (i)} $M_i$,  $i=1,\ldots,m-1$, is the monodromy matrix of a local system (\ref{si}) at $a_i$ realizing the local Galois data $(G_i,\ll_i)$ with true Poincar\'e rank $r_i$,

\noindent {\it (ii)} the representation $\chi_{\mb}$  is irreducible,

\noindent then $G$ is the Galois group over $\cc(z)$ of a ssystem (\ref{gs}) with singularities all in $\cd$, whose Poincar\'e rank  at each $a_i$ is $r_i$,  $i=1,\ldots,m-1$,  and which is fuchsian at~$a_m$ . 
\end{corollary}

\medskip
\noindent If a  group with a local Galois structure is  {\it connected},  then it has  local structures of the form $\ll=({\mathrm {id}}, T, \cn)$.  It is in then possible (\cf \cite{R4}, 4.2.2)  to realize $\ll$ with a system for which the formal and the topological monodromy matrices  coincide and are  sequal to ${\mathrm {exp}}(u)$,   for any generator $u$ of the Lie algebra $\cn$. We obtain the following result in this case.

\begin{corollary} \label{corgaloisrecol}  Let $G$ be a linear algebraic group over $\cc$. Assume that $G$ is topologically generated by closed connected subgroups  $G_1,\ldots,G_{m-1}$, $m\ge 2$, each endowed with a reduced local Galois structure $\ll_i=({\mathrm {id}}, T_i, \cn_i)$  and let $\cd=\{a_1,\ldots,a_m\}$  be  an arbitrary set of $m$  points of $\pp$.  Let $u_i$, for $i=1,\ldots,m-1$, be a generator of the Lie algebra $\cn_i$, and $r_i$   the true Poincar\'e of a local system (\ref{si}) realizing  $(G_i,\ll_i)$ with the monodromy matrix $M_i={\mathrm {exp}}(u_i)$. If the representation $\chi_{\mb}$, where  $\mb=(M_1,\ldots,M_{m-1})$, is irreducible, then $G$ is the Galois group over $\cc(z)$ of a system (\ref{gs}) with singularities all in $\cd$, whose Poincar\'e rank  at each $a_i$ is $r_i$, $i=1,\ldots,m-1$,  and  which is fuchsian at $a_m$. 
\end{corollary}

\noindent
{\bf Remark} : Corollaries \ref{galoisirr} and \ref{corgaloisrecol} in particular hold if, for some or all $i$,  the local systems in the statements  realize the {\it minimal} Poincar\'e rank $r_{{\ll}_i}$ at $a_i$.

\medskip
\noindent We  now wish to realize a given group globally with $\sb(G)$ singularities, all fuchsian but possibly one,  with a minimal Poincar\'e rank at the irregular singularity.

\medskip
\noindent
We will  use yet another characterization of a local Galois group $\Gamma$ (\cf  \cite{PS},  Theorem 11.13), namely that $V(\Gamma)=\Gamma/L(\Gamma)$ be topologically generated by one element.

\medskip
\noindent  It follows  from this criterion that for any $\alpha\in V(G)$ the inverse image $G_{\alpha}=pr^{-1}(<\alpha>)$, by the projection $pr : G\rightarrow V(G)$,  of the closed subgroup topologically generated by $\alpha$, has a local Galois structure.

\medskip
\noindent
{\bf Notation} {\it Let $\aa$ denote the set of all elements $\alpha\in V(G)$ which belong to a minimal family of generators of $V(G)$, and let $r(G)={\mathrm{min}}_{\alpha\in \aa}(r(G_{\alpha}))$ denote the minimal possible Poincar\'e rank  of a  system realizing  local Galois data $(G_{\alpha},\ll_{\alpha})$, $\alpha\in \aa$. 

\noindent Let us  write $\sb$ for $\sb(G)$, and $\overline{X}$ for the class of an element $X$ of $G$ in the quotient $V(G)$. }

\medskip
\noindent With this notation,  we obtain the following result.
\begin{theorem}  \label{irrmin} Let $\cd=\{a_1,\ldots,a_{\sb}\}$  be an arbitrary set of $\sb$  points of $\pp$. Let $M_1$ denote the monodromy matrix  of a system of  Poincar\'e rank $r(G)$ at $a_1$  and realizing $G_{\alpha}$, for some $\alpha\in \aa$, as its local Galois group. Consider a subset  $\{M_2,\ldots,M_{\sb -1}\}\subset G$ such that $(\alpha,\overline{M_2},\ldots,\overline{M_{\sb-1}})$ is a minimal family of generators of $V(G)$. Let $\mm$ be a reduced datum on $\cd$ which includes the representation $\chi_{\mb}$ for $\mb=(M_1, M_2,\ldots,M_{{\sb}-1})$,  fuchsian data at $a_2,\ldots,a_{\sb},$ and  the Poincar\'e rank $r(G)$  at $a_1$. If $\mm$ fulfills the conditions of Theorem \ref{main}, then $G$ is the  Galois group over $\cc(z)$ of a system with fuchsian singularities at $a_2,\ldots,a_{\sb}$ and  Poincar\'e rank $r(G)$ at $a_1$.
\end{theorem}

\begin{proof} The group $G$ is  generated by the closed subgroups $G_{\alpha}$ and $G_2,\ldots,G_{\sb-1}$, where $G_i$, $i=2,\ldots,\sb-1$, is  generated by $M_i$.  It is then possible to define generalized monodromy data  which include $\chi_{\mb}$ and a true Poincar\'e rank equal to $r(G)$ at $a_1$. Apply Theorem \ref{galoisrecol} to conclude. 
\end{proof}

\medskip
\noindent
The irreducibility  condition on the monodromy in particular implies the following result.

\begin{corollary} 
If for some $\alpha\in\aa$  such that $r(G_{\alpha})=r(G)$ and for a family
$\mb=(M_2,\ldots,M_{\sb -1})$  of elements of $G$ such that $(\alpha,\overline{M_2},\ldots,\overline{M_{\sb-1}})$ is a minimal family of generators of $V(G)$ the representation $\chi_{\mb}$ is irreducible, then $G$ is the  Galois group over $\cc(z)$ of a system with no more than $\sb(G)$ singularities, all fuchsian but  one, possibly irregular, at which the system has Poincar\'e rank  $r(G)$.
\end{corollary}

\begin{proof} The group $G$ is  generated by the closed subgroup $G_{\alpha}$ and $G_2,\ldots,G_{\sb-1}$, where $G_i$, $i=1,\ldots,\sb-1$, is  generated by $M_i$. Consider a subset $\cd=\{ a_1,\ldots,a_{\sb} \}$ of $\pp$ and let $M_1$ denote the monodromy matrix  of a system of minimal Poincar\'e rank $r(G)$ at $a_1$  realizing $G_{\alpha}$ as its local Galois group. The representation  $\chi_{\mb'}$, where $\mb'=(M_1,\ldots,M_{\sb-1})$, is irreducible, and we can apply Corollary \ref{galoisirr} to conclude. 
\end{proof}

\noindent We can restate this result with a weaker condition on the monodromy representation, but with less control on the Poincar\'e rank.
\begin{corollary}  Consider a subset  $\cd=\{ a_1,\ldots,a_{\sb} \}$ of $\pp$. Assume that for some $\alpha\in\aa$ and a family $\mb=(M_1,\ldots,M_{\sb -1})$ of elements of $G$ one has

\smallskip
\noindent
{\it (i)} $\overline{M_1}=\alpha$,

\smallskip
\noindent
{\it (ii)} $(\overline{M_1},\ldots,\overline{M_{\sb-1}})$ is a minimal family of generators of $V(G)$,

\smallskip
\noindent
{\it (iii)} the representation $\chi_{\mb}$ is irreducible.

\smallskip
\noindent Let $r_1$ be the true Poincar\'e rank at $a_1$ of a system  realizing $G_{\alpha}$ as its local Galois group at $a_1$ and with   monodromy matrix $M_1$. Then $G$ is the  Galois group over $\cc(z)$ of a system with no more than $\sb(G)$ singularities, all fuchsian but  one, possibly irregular, at which the system has Poincar\'e rank  $r_1$.

\end{corollary}
\begin{proof}
The existence of  a system (\ref{si}) with the monodromy matrix $M_1\in G_{\alpha}$  at $a_1$ follows from Ramis's construction for the solution of the local inverse problem  (\cf \cite{PS}, section 11).  We may assume that the Poincar\'e rank of (\ref{si}) at $a_1$ is minimal (equal to the true Poincar\'e rank). This local system at $a_1$, together with the representation $\chi_{\mb}$ and fuchsian data at $a_2,\ldots,a_{\sb}$, defines a reduced datum $\mm$ which by Corollary \ref{irr} can be realized with a global system (\ref{gs}) since $\chi_{\mb}$ is irreducible and one point at least, $a_{\sb}$, is without roots. The Galois group of  (\ref{gs}) is clearly $G$, and its Poincar\'e rank at $a_1$ is $r_1$ by the generalized Riemann-Hilbert problem.
\end{proof}
\bigskip
\noindent
In dimension two and three we can say more. To apply the  results for the GRH-problem to this case,  we  need to assume  that the generalized monodromy data are, at the irregular singularity, those of a local system with a divergent fundamental solution. 

\medskip
\noindent
We  first recall a characterization of  Galois groups over the differential field $\cc((z))$ of formal Laurent series, or {\it formal Galois groups}. Ramis has characterized such groups by a {\it formal local Galois structure} $(T,a,\cn)$ which only differs from  the above mentionned  local Galois structure by the  condition {\it (iii)}, which in the formal case is replaced by  $$\gg=\cn + \ct \leqno{(iii)'}$$
where $\ct$ denotes the Lie algebra of the torus $T$. In (\cite{PS}, Theorem 11.2) formal local Galois groups $\Gamma$ are characterized by the  simpler, equivalent condition that $\Gamma$ contain a normal subgroup $T$ such that $T$ is a torus and $\Gamma/T$ be topologically generated by one element.

\medskip
\noindent
If we apply Theorem \ref{dt} of section 5  in dimension two and three we obtain the following result.
\begin{theorem}
Let $G$ be a linear algebraic subgroup of $\gl(p,\cc)$, $p=2, 3$. If

\noindent (i) \  p=2, or

\noindent (ii)  p=3 and  for some $\alpha\in \aa$  such that $r(G_{\alpha})=r(G)$, $G_{\alpha}$  is not a formal Galois group, 

\noindent then $G$ is the Galois group over $\cc(z)$ of a linear differential system with no more than $\sb(G)$ singularities, all fuchsian but one,  irregular of Poincar\'e rank $r(G)$.
\end{theorem}

\begin{proof} Since, if $p=3$, the subgroup $G_{\alpha}$ is not a formal Galois group, any of its local Galois structures  will realize it as the Galois group (over $\cc(\{z\})$) of a local system with a divergent fundamental solution. Otherwise, the Stokes matrices would be trivial, hence the formal solutions would be convergent. Moreover, since $\sb(G)\ge 2$ we can assume, for $p=2$ and $3$, that the data at the other singularities are all fuchsian, and apply the results of Corollaries 4 and 5  to conclude. \end{proof}

\bigskip
\noindent
{\bf Note} Andrey Andreevich Bolibruch died during the completion of this paper. We dedicate it to the memory of our
late coauthor, colleague, and wonderful friend.

\end{document}